\documentclass[10pt]{amsart}
\usepackage{amscd,amssymb}

\catcode`@=11 \catcode`!=11

\expandafter\ifx\csname fiverm\endcsname\relax
\expandafter\ifx\csname fivrm\endcsname\relax
   \relax
\else
  
\fi
\fi

\let\!latexendpicture=\endpicture 
\let\!latexframe=\frame
\let\!latexlinethickness=\linethickness
\let\!latexmultiput=\multiput
\let\!latexput=\put
 
\def\@picture(#1,#2)(#3,#4){%
  \@picht #2\unitlength
  \setbox\@picbox\hbox to #1\unitlength\bgroup 
  \let\endpicture=\!latexendpicture
  \let\frame=\!latexframe
  \let\linethickness=\!latexlinethickness
  \let\multiput=\!latexmultiput
  \let\put=\!latexput
  \hskip -#3\unitlength \lower #4\unitlength \hbox\bgroup}

\catcode`@=12 \catcode`!=12

\input pictex
\def\begingraph<#1,#2>[#3]{%
  \beginpicture
  \setcoordinatesystem units <#1, #2>
  \setlinear
  \linethickness=#3
  }
\def\endgraph{%
  \endpicture
  }
\def\beginsmallgraph{\begingraph<0.5cm,0.5cm>[0.35pt]}

\def\vertex at (#1, #2){\put{$\bullet$} at #1 #2}
\catcode`!=11 
\def\edge from (#1, #2) to (#3, #4){%
  \!xloc=\!M{#1}\!xunit  \!xxloc=\!M{#3}\!xunit%
  \!yloc=\!M{#2}\!yunit  \!yyloc=\!M{#4}\!yunit%
  \!dxpos=\!xxloc  \advance\!dxpos by -\!xloc
  \!dypos=\!yyloc  \advance\!dypos by -\!yloc
  \ifdim\!dypos=\!zpt
    \putrule from #1 #2 to #3 #4
  \else
    \ifdim\!dxpos=\!zpt
      \putrule from #1 #2 to #3 #4
    \else
      \plot #1 #2 #3 #4 /
    \fi
  \fi
  }
\catcode`!=12 
\def\Label #1 at (#2, #3){\put {#1} at #2 #3}
\catcode`@=11 \catcode`!=11
  
\let\!pictexendpicture=\endpicture 
\let\!pictexframe=\frame
\let\!pictexlinethickness=\linethickness
\let\!pictexmultiput=\multiput
\let\!pictexput=\put

\def\beginpicture{%
  \setbox\!picbox=\hbox\bgroup%
  \let\endpicture=\!pictexendpicture
  \let\frame=\!pictexframe
  \let\linethickness=\!pictexlinethickness
  \let\multiput=\!pictexmultiput
  \let\put=\!pictexput
  \let\input=\@@input   
  \!xleft=\maxdimen  
  \!xright=-\maxdimen
  \!ybot=\maxdimen
  \!ytop=-\maxdimen}

\let\frame=\!latexframe

\let\pictexframe=\!pictexframe

\let\linethickness=\!latexlinethickness
\let\pictexlinethickness=\!pictexlinethickness

\let\\=\@normalcr
\catcode`@=12 \catcode`!=12

\theoremstyle{plain}
\newtheorem{thm}[subsection]{Theorem}

\newtheorem{prop}[subsection]{Proposition} 
\newtheorem{cor}[subsection]{Corollary}

\theoremstyle{definition}
\newtheorem{rem}[subsection]{Remark}
\newtheorem{defn}[subsection]{Definition} 
\newtheorem{exmp}[subsection]{Example}

\numberwithin{equation}{section}

\newenvironment{alphenum}%
{%
 \begin{enumerate}
}
{%
 \end{enumerate}%
}

\newenvironment{romenum}%
{%
 \begin{enumerate}
}
{%
 \end{enumerate}%
}

\newcommand{\A}{\mathcal{A}}
\newcommand{\B}{\mathcal{B}}
\newcommand{\V}{\mathcal{V}}
\newcommand{\CC}{\mathcal{C}}
\newcommand{\LL}{\mathcal{L}}
\newcommand{\KK}{\mathcal{K}}
\newcommand{\MM}{\mathcal{M}}
\newcommand{\Az}{\mathcal{A_{\omega}}}
\newcommand{\C}{\mathbb{C}}
\newcommand{\F}{\mathbb{F}}

\newcommand{\Z}{\mathbb{Z}}
\newcommand{\R}{\mathbb{R}}
\newcommand{\RP}{\mathbb{RP}}
\newcommand{\SP}{\mathbb{S}}
\newcommand{\D}{\Delta}
\newcommand{\e}{\epsilon}

\newcommand{\f}{\phi}

\newcommand{\s}{\sigma}

\newcommand{\T}{\Theta}

\newcommand{\bt}{\mathbf{t}}

\renewcommand{\a}{{\alpha }}
\renewcommand{\b}{{\beta }}

\renewcommand{\l}{\lambda}
\renewcommand{\t}{\tau}

\renewcommand{\L}{{\Lambda }}

\DeclareMathOperator{\id}{id}

\DeclareMathOperator{\Inn}{Inn}
\DeclareMathOperator{\Center}{Center}
\DeclareMathOperator{\GL}{GL}
\DeclareMathOperator{\ab}{ab}

\DeclareMathOperator{\lk}{lk}
\def\Re{\operatorname{Re}}
\def\Im{\operatorname{Im}}
\DeclareMathOperator{\Tors}{Tors}
\DeclareMathOperator{\sgn}{sgn}
\DeclareMathOperator{\depth}{depth}

\DeclareMathOperator{\ii}{i}

\begin{document}

\title[Homotopy Types of $2$-Arrangements]%
{Homotopy Types of Complements of \\
$2$-Arrangements in $\R^4$}

\author[D.~Matei]{Daniel~Matei}
\address{Department of Mathematics,
Northeastern University,
Boston, MA 02115}
\email{dmatei@lynx.neu.edu}

\author[A.~Suciu]{Alexander~I.~Suciu$^{\dag}$}
\address{Department of Mathematics,
Northeastern University,
Boston, MA 02115}
\email{alexsuciu@neu.edu}
\urladdr{http:\slash\slash www.math.neu.edu\slash\~{}suciu}

\thanks{$\dag$Partially supported 
by N.S.F.~grant DMS--9504833.}

\subjclass{Primary 57M05, 57M25, 52B30; Secondary 14M12, 20F36} 

\keywords{arrangement, line configuration, link, 
braid, characteristic variety}

\begin{abstract}
We study the homotopy types of complements of arrangements 
of $n$ transverse planes in $\R^4$, obtaining a complete 
classification for $n\le 6$, and lower bounds for the number 
of homotopy types in general.  Furthermore, we show that the 
homotopy type of a $2$-arrangement in $\R^4$ is not determined 
by the cohomology ring, thereby answering a question of Ziegler. 
The invariants that we use are derived from the characteristic 
varieties of the complement.  The nature of these varieties 
illustrates the difference between real and complex arrangements.  
\end{abstract}

\maketitle


\section{Introduction}  
\label{sec:intro}

In~\cite{GM}, Goresky and MacPherson introduced a generalization of the
notion of complex hyperplane arrangement. A {\em $2$-arrangement} 
in $\R^{2d}$ is a finite collection $\A$ of codimension~$2$ linear 
subspaces so that, for every subset $\B\subseteq \A$, the space 
$\bigcap_{H\in \B}\, H$ has even dimension. 
The main object of study is the
complement of the arrangement, $X(\A)=\R^{2d}\setminus \bigcup_{H\in \A}\, H$. 
Goresky and MacPherson computed the cohomology groups of $X$.  Bj\"{o}rner
and Ziegler \cite{BZ} and Ziegler \cite{Z} determined the structure of the
cohomology algebra $H^*(X;\Z)$. These results generalize the classical work
of Arnol'd, Brieskorn, and Orlik and Solomon on the cohomology ring of the
complement of a complex hyperplane arrangement, see \cite{OT}.  
Unlike the situation obtaining for the Orlik-Solomon algebra, which is 
completely determined by the intersection lattice, 
there remained an ambiguity in the relations defining $H^*(X;\Z)$. 
Even in the simplest case of $2$-arrangements in $\R^4$,  
a striking phenomenon occurs, showing that this ambiguity 
cannot be resolved,~\cite{Z}.  

Each $2$-arrangement $\A$ in $\R^4$ is a realization of the uniform 
matroid $U_{2,n}$, where $n=|\A|$ is the cardinality of the arrangement.  
Thus, the intersection lattice of such an arrangement is uniquely 
determined by $n$.  Furthermore, the homology groups of the complement, 
$X$, the lower central series quotients of the group $G=\pi_1(X)$, 
and the Chen groups of $G$ also depend only on $n$.

On the other hand, the cohomology ring of $X$ is a more subtle invariant.  
The relations in $H^*(X;\Z)$ depend on, and are determined by the 
linking numbers of the associated link.  Ziegler ~\cite{Z} found a pair 
of $2$-arrangements of four planes which have non-isomorphic 
cohomology rings. His method, which uses an invariant derived 
from $H^*(X;\Z)$, does not seem, however, to extend beyond $n=4$.  

In this paper, we introduce new homotopy-type invariants 
of complements of $2$-arrangements.  These invariants, derived 
from the Alexander module, work for arbitrary $n$.  As a first step 
towards the homotopy classification of $2$-arrangements, we prove 
the following (see Corollary~\ref{cor:depth2}).  

\begin{thm}  \label{thm:lowerbound}
For every integer $n\ge 1$, there exist at least 
$p(n-1) - \lfloor\frac{n-1}{2}\rfloor$ 
different homotopy types of complements of 
$2$-arrangements of $n$~planes in $\R^4$, 
where $p( \cdot)$ is the partition 
function, and $\lfloor\cdot\rfloor$ 
is the integer part function.    
\end{thm}

At the end of \cite{Z}, Ziegler asks whether the cohomology 
ring determines the homotopy type of the complement of a $2$-arrangement, 
proposing as a candidate for a negative answer the remarkable pair 
of arrangements of $6$ planes found by Mazurovski\u{\i} in \cite{M1}.  
Using successive cablings on Mazurovski\u{\i}'s pair, we answer Ziegler's 
question, as follows (see Theorem~\ref{thm:mazcables}). 

\begin{thm}  \label{thm:homotopytypes}
For every integer $n\ge 6$, there exists a pair of $2$-arrangements 
of $n$ planes in $\R^4$, whose complements have isomorphic cohomology rings, 
but different homotopy types.  
\end{thm}

Rigid isotopy of arrangements implies isotopy of their singularity 
links.  The converse is not clear, though, since an isotopy may go outside 
the class of such links.  On the other hand, the classification, 
up to rigid isotopy, of $2$-arrangements in $\R^4$  is equivalent 
to the classification, also up to rigid isotopy, of configurations 
of skew lines in $\R^3$.  Such configurations were introduced by 
Viro in \cite{V}, and have been intensively studied  since then, 
see the survey article by Crapo and Penne~\cite{CP}.  The rigid 
isotopy classification of configurations of $n$ skew lines in $\R^3$, 
was achieved by Viro~\cite{V} for $n\le 5$, and by 
Mazurovski\u{\i}~\cite{M1} for $n=6$.  

It is readily seen that rigid isotopy of arrangements implies homotopy 
equivalence of their complements.  The converse is not true.  Indeed, 
as first noted by Viro, there exist configurations that 
are not isotopic to their mirror image.  But clearly, the complements 
of mirror pairs are diffeomorphic, and thus homotopy equivalent.  
The next result shows that this is the only exception, for $n\le 6$ 
(see Theorem~\ref{thm:classify6}).

\begin{thm}  \label{thm:lowclassification}  
For $2$-arrangements of $n\le 6$ planes in $\R^4$, 
the homotopy types of complements are in one-to-one 
correspondence with the rigid isotopy types modulo mirror images.  
\end{thm}

This theorem recovers Ziegler's classification 
of homotopy types of arrangements of $n=4$ planes.
The number of homotopy types from the classification in 
Theorem~\ref{thm:lowclassification}, together with the lower bound 
from Theorem~\ref{thm:lowerbound}, are tabulated 
below.\footnote[2]{Recently, Borobia and Mazurovski\u{\i}~\cite{BM} achieved 
the rigid isotopy classification of configurations of $7$ lines.  
If the assertion of Theorem~\ref{thm:lowclassification} were 
to hold for $n=7$, it would give $37$ distinct homotopy types 
of complements of arrangements of $7$ planes.  We have verified 
this in the particular case of horizontal arrangements, for which there 
are $24$ distinct homotopy types.}

\[
\begin{array}{|l||c|c|c|c|c|c|c|}
\hline
\multicolumn{1}{|c||}{n}
&\multicolumn{1}{|c|}{1}
&\multicolumn{1}{|c|}{2}
&\multicolumn{1}{|c|}{3}
&\multicolumn{1}{|c|}{4}
&\multicolumn{1}{|c|}{5}
&\multicolumn{1}{|c|}{6}
&\multicolumn{1}{|c|}{7}
\\ 
\hline\hline
\text{Homotopy types}
& 1 & 1 & 1 & 2 & 4 & 11 & ?\\
\hline
\text{Lower bound} 
& 1 & 1 & 1 & 2 & 3 & 5 & 8 \\
\hline
\end{array}
\]

The homotopy-type invariants that we use in our classification 
of $2$-arrangements are derived from the characteristic varieties of 
their complements.  Given a space $X$ with $H_1(X)\cong \Z^n$, 
the $k^{\text{th}}$ determinantal ideal of the Alexander module 
of $X$ defines a subvariety, $V_k(X)$, of the complex algebraic 
torus, $(\C^*)^n$, whose monomial isomorphism type depends 
only on the homotopy type of $X$---in fact, only on 
$\pi_1(X)$---see~\cite{DF, L1}.  
We call $V_k(X)$ the $k^{\text{th}}$ characteristic variety of $X$.  
From this variety, we extract in Theorem~\ref{thm:invariant} 
the following homotopy-type invariants for the space $X$: 
the list $\Sigma_k(X)$ of codimensions of irreducible components,  
and the number $\Tors_{p,k}(X)$ of $p$-torsion points. 
These numerical invariants are readily computable by standard 
methods of geometric topology and commutative algebra, and are 
powerful enough to detect all the differences in homotopy 
types listed in the above theorems.  

The characteristic varieties of complements of divisors in complex 
algebraic manifolds have been intensively studied recently, 
see~\cite{Ar, L2, Hr, CScv, L3, LY}.  
Deep results as to their qualitative nature have 
been obtained by Arapura \cite{Ar}, who showed that all the 
irreducible components of such characteristic varieties are 
(possibly translated) subtori of a complex algebraic torus.  
Building on this work, a more precise description of the 
characteristic varieties of complex hyperplane arrangements  
has emerged.  In all known examples, if $X$ is the complement 
of such an arrangement, all positive-dimensional subtori of 
$V_k(X)$ pass through the origin $\mathbf{1}$ of the torus.  

On the other hand, if $X$ is the complement of a $2$-arrangement  
in $\R^4$, we find that the characteristic varieties of $X$ may 
contain positive-dimensional subtori that do not pass through~$\mathbf{1}$.  
For the non-complex Ziegler arrangement, the variety $V_2$ contains three 
subtori of $(\C^*)^4$, one of which is translated by $(1,-1,1,1)$, 
see Example~\ref{ex:Zieglercharvar}.  But this is still a rather mild 
qualitative difference.  For the indecomposable Mazurovski\u{\i} 
arrangements, the variety $V_1$ is not even a union of translated 
subtori, see \S\ref{sec:maz}.  These phenomena may be thought 
of as manifestations of the non-complex nature of real arrangements.  

The paper is organized as follows.  

In \S\ref{sec:arrconlin}, we review the basic facts about 
$2$-arrangements in $\R^4$, and their associated configurations of 
lines and singularity links.  
In \S\ref{sec:decomphoriz}, 
we look in detail at some special classes of arrangements: 
the decomposable ones, and the horizontal ones.  
In \S\ref{sec:braids}, we associate several braids to a 
$2$-arrangement, and use these braids to compute the fundamental 
group of the complement.  
In \S\ref{sec:detideals}, we review Alexander modules and 
define numerical homotopy-type invariants from the associated 
characteristic varieties.  
In \S\ref{sec:bottomvar}, we study the bottom 
characteristic varieties $V_{n-2}$, obtaining a complete characterization 
for depth~$2$, completely decomposable arrangements.  
In \S\ref{sec:topvar}, we study the top characteristic varieties $V_1$, 
and their torsion points. 
In \S\ref{sec:maz}, we study in detail the Mazurovski\u{\i} 
arrangements, and their cablings.  
Using the results and techniques from 
\S\S\ref{sec:bottomvar}--\ref{sec:maz}, we complete the homotopy-type 
classification of $2$-arrangements of $6$ planes or less in 
\S\ref{sec:classification}. 

\medskip
\noindent
{\bfseries Acknowledgment.}   This work started from an illuminating 
conversation with G\"{u}nter Ziegler, who introduced us to~\cite{Z}.  
The computational part of the work was greatly aided by 
{\it Mathematica${^{\circledR}}$}, and by the commutative 
algebra package {\it Macaulay~2}.  Thanks are due to the referee, 
for many valuable suggestions that have improved both the substance 
and the style of the paper.


\section{Arrangements, Line Configurations, and Links} 
\label{sec:arrconlin} 

In this section we collect some facts about arrangements of transverse 
planes in $\R^4$, and the corresponding configurations of skew lines 
in $\R^3$ and links in $\SP^3$.  

\subsection{}  \label{defpoly}
We start by defining our basic objects of study in a concrete way.  

\begin{defn}  \label{def:2arr}
A {\em $2$-arrangement} in $\R^4$ is a finite collection 
$\A=\{H_1,\dots ,H_n\}$ of pairwise transverse $2$-dimensional 
vector subspaces of $\R^{4}$.  The {\em union} of the arrangement 
is $U(\A)=\bigcup_{H\in \A} H$.  The {\em complement} of the arrangement 
is $X(\A)=\R^4\setminus U(\A)$.  The {\em link} of the arrangement is 
$L(\A)=\SP^{3}\cap U(\A)$. 
\end{defn}

Each plane $H_i\in \A$ can be written as 
$H_i=\ker\lambda_i \cap \ker\lambda'_i$, for some linear forms 
$\lambda_i, \lambda'_i:\R^4\to \R$.  
The transversality condition means that $H_i\cap H_j=\{0\}$, 
for all $i\ne j$.  That is, 
$\det(\l_i,\l'_i,\l_j,\l'_j)\ne 0$ for $i\ne j$. 

Alternatively, identifying $\R^4$ with $\C^2=\{(z,w)\}$, 
each plane in $\A$ can be written as $H_i= \{f_i=0\}$, where 
$f_i(z,w)=a_iz+b_i\bar{z}+c_iw+d_i\bar{w}$, 
for some $a_i,b_i,c_i,d_i\in\C$.  
In terms of real coordinates 
$x=\Re z$, $y=\Im z$, $u=\Re w$, $v=\Im w$, we have 
$\lambda_i(x,y,u,v)=\Re f_i(x+\ii y, u+\ii v)$ and 
$\lambda'_i(x,y,u,v)=\Im f_i(x+\ii y, u+\ii v)$, where $\ii=\sqrt{-1}$.  

With notation as above, let $f:\C^2\to \C$ be the polynomial 
map in $z, \bar{z}, w, \bar{w}$ given by $f=f_1\cdots f_n$.  We say 
that $f$ is a {\em defining polynomial} for the arrangement $\A$.  
Obviously, the union of the arrangement is the zero locus of the 
defining polynomial. 

\begin{exmp}  \label{complexlines}
The most basic example of a $2$-arrangement is a {\em complex} arrangement. 
Such an arrangement consists of complex lines through the origin of $\C^2$.  
 Any two complex arrangements differ by an 
$\R$-linear change of variables, and thus have diffeomorphic complements.  
We denote the complex arrangement of $n$ lines by $\A_n$, 
and take its defining polynomial to be $f_n(z,w)=(z-w)\cdots (z-nw)$.  
The link $L(\A_n)$ is the $n$-component Hopf link.  
The {\em trivial} arrangement is $\A_1$. 
\end{exmp}

\subsection{}  \label{singlink}
Let $\A=\{H_1,\dots ,H_n\}$ be a $2$-arrangement in $\R^4$.  
Its link, $L=\{L_1,\dots , L_n\}$, consists of $n$ unknotted circles 
in $\SP^3$.  The complement of the arrangement, $X(\A)$, is homotopy 
equivalent to the complement of the link, $Y(L)=\SP^{3}\setminus L$, 
via radial deformation.  Using this observation, we can compute 
homotopy-type invariants of $X=X(\A)$ by methods of knot theory.  

The homology groups of $X$ depend only on the number 
of planes in the arrangement:  $H_0=\Z$, $H_1=\Z^n$, $H_2=\Z^{n-1}$, 
$H_k=0$ for $k>2$.  The cohomology ring of $X$, on the other hand, 
also depends on the linking numbers $l_{i,j}=\lk(L_i,L_j)$.  
Specifically, 
\begin{equation*}  \label{eq:cohomology}
H^{*}(X; \Z) = \left.\sideset{}{^*}\bigwedge\Z^n \right\slash \left( 
l_{ij}e_ie_j + l_{jk}e_je_k + l_{ki}e_ke_i = 0 \right),
\end{equation*}
where $\bigwedge^{*} \Z^n$ is the exterior algebra on $e_1, \dots , e_n$. 
As noted by Ziegler~\cite{Z}, one can compute the linking numbers of 
$L(\A)$ directly from the defining equations of $\A$.  Indeed, if 
$H_i=\{\lambda_i=\lambda'_i=0\}$, then 
$l_{i,j}=\sgn(\det(\lambda_i,\lambda'_i,\lambda_j,\lambda'_j))$.  

As shown by Ziegler~\cite{Z}, the complement $X$ fibers over 
$\C^*=\C\setminus \{0\}$, with fiber $\C\setminus \{n-1\ \text{points}\}$, 
and thus $X$ is a $K(G,1)$ space.  Alternatively, since all the linking 
numbers are non-zero, the link $L$ is non-split, and thus $Y(L)$ 
is aspherical, see~\cite{BuZ}.  It follows that the homotopy type of $X$ 
is determined by the isomorphism class of its fundamental group $G$. 

As we shall see in Proposition~\ref{prop:bundle}, the monodromy 
of the bundle $X\to \C^*$ is a certain (pure) braid automorphism 
$\check{\beta}\in P_{n-1}$, and so $G$ is a semidirect product of 
free groups, $G=\F_{n-1}\rtimes_{\check{\beta}} \F_1$.  Since 
$\check{\beta}$ acts trivially on homology, a result of Falk and 
Randell~\cite{FR} implies that the lower central series quotients 
of $G$ depend only on $n$, being equal to those of the product 
$\Gamma=\F_{n-1}\times \F_1$.  In fact, since all the linking 
numbers of $L$ are equal to $\pm 1$, a result of Massey and Traldi~
~\cite{MT} shows that the lower central series quotients of both 
$G$ and $G/G''$ are equal to the corresponding quotients of 
$\Gamma$ and $\Gamma/\Gamma''$. 

\subsection{}  \label{config} 
Now let $H$ be an affine hyperplane in $\R^4$, generic with respect 
to $\A$.  The {\em configuration} of $\A$ corresponding to $H$ is 
the configuration of skew lines in $\R^{3}$ defined as 
$\CC_{H}(\A)=H\cap U(\A)$.  Conversely, given a configuration $\CC$ 
of skew lines, one obtains a $2$-arrangement, $\A_p=\A(\CC)$, by 
coning at a generic point $p$, and translating $p$ to $0$.

\begin{figure}
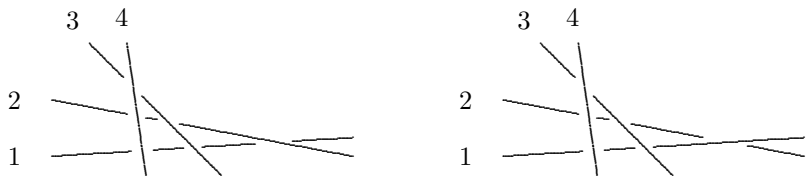

\centerline{\hfill\hskip -1.28 in
\beginsmallgraph
\Label {$1$} at (1, 1)
\Label {$2$} at (1, 2.5)
\Label {$3$} at (2.55, 4.6)
\Label {$4$} at (3.85, 4.7)
\edge from (4, 4) to (4.5, 0.5)
\edge from (3, 4) to (3.9, 3.1)
\edge from (4.4, 2.6) to (6.5, 0.5)
\edge from (2, 2.5) to (4, 2.125)
\edge from (4.5, 2.03125) to (4.8, 1.975)
\edge from (5.4, 1.8625) to (10, 1)
\edge from (2, 1) to (4.1, 1.13125)
\edge from (4.7, 1.16875) to (5.5, 1.21875)
\edge from (6, 1.25) to (7.4, 1.3375)
\edge from (8.4, 1.4) to (10, 1.5)
\Label {$1$} at (13, 1)
\Label {$2$} at (13, 2.5)
\Label {$3$} at (14.55, 4.6)
\Label {$4$} at (15.85, 4.7)
\edge from (16, 4) to (16.5, 0.5)
\edge from (15, 4) to (15.9, 3.1)
\edge from (16.4, 2.6) to (18.5, 0.5)
\edge from (14, 2.5) to (16, 2.125)
\edge from (16.5, 2.03125) to (16.8, 1.975)
\edge from (17.4, 1.8625) to (19.3, 1.50625)
\edge from (20.5, 1.28125) to (22, 1)
\edge from (14, 1) to (16.1, 1.13125)
\edge from (16.7, 1.16875) to (17.5, 1.21875)
\edge from (18, 1.25) to (22, 1.5)
\endgraph
\hfill}
\caption{Ziegler's pair: 
The configurations $\CC^{+}$ (left) and $\CC^{-}$ (right).}
\label{fig:zpair}
\end{figure}

\begin{exmp} \label{ex:ziegler} 
Let $\A^{+}$ and $\A^{-}$ be the pair of $2$-arrangements considered 
by Ziegler in~\cite{Z}.  The arrangement $\A^{+}$ is the complex 
arrangement $\A_4$, with defining polynomial 
$f^{+}(z,w)=(z-w)(z-2w)(z-3w)(z-4w)$. 
The arrangement $\A^{-}$ has defining polynomial 
$f^{-}(z,w)=(z-\bar{w})(z-2\bar{w})(z-3w)(z-4w)$. 
Projecting onto the hyperplane $\{v=1\}$, we get configurations 
$\CC^{\pm}=\{\ell_1^{\pm},\ell_2^{\pm},\ell_3,\ell_4\}$, 
with equations
\begin{alignat*}{2}  \label{eq:4lines}
\ell_1^{\pm}&=\{x-u=y\mp 1=0\}, 
&\qquad \ell_2^{\pm}&=\{x-2u=y\mp 2=0\},\\ 
\ell_3&=\{x-3u=y-3=0\},
&\qquad \ell_4&=\{x-4u=y-4=0\}. 
\end{alignat*}
\end{exmp}
\noindent
The two configurations are pictured in Figure~\ref{fig:zpair}.

\begin{figure}
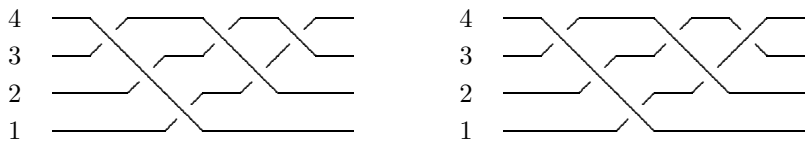

\centerline{\hfill  \hskip -0.1 in
\beginsmallgraph
\Label {$1$} at (1, 1)
\Label {$2$} at (1, 2)
\Label {$3$} at (1, 3)
\Label {$4$} at (1, 4)
\edge from (2, 1) to (3, 1)
\edge from (2, 2) to (3, 2)
\edge from (2, 3) to (3, 3)
\edge from (2, 4) to (3, 4)
\edge from (3, 1) to (4, 1)
\edge from (3, 2) to (4, 2)
\edge from (3, 3) to (3.3, 3.3)
\edge from (3.7, 3.7) to (4, 4)
\edge from (3, 4) to (4, 3)
\edge from (4, 1) to (5, 1)
\edge from (4, 2) to (4.3, 2.3)
\edge from (4.7, 2.7) to (5, 3)
\edge from (4, 3) to (5, 2)
\edge from (4, 4) to (5, 4)
\edge from (5, 1) to (5.3, 1.3)
\edge from (5.7, 1.7) to (6, 2)
\edge from (5, 2) to (6, 1)
\edge from (5, 3) to (6, 3)
\edge from (5, 4) to (6, 4)
\edge from (6, 1) to (7, 1)
\edge from (6, 2) to (7, 2)
\edge from (6, 3) to (6.3, 3.3)
\edge from (6.7, 3.7) to (7, 4)
\edge from (6, 4) to (7, 3)
\edge from (7, 1) to (8, 1)
\edge from (7, 2) to (7.3, 2.3)
\edge from (7.7, 2.7) to (8, 3)
\edge from (7, 3) to (8, 2)
\edge from (7, 4) to (8, 4)
\edge from (8, 1) to (9, 1)
\edge from (8, 2) to (9, 2)
\edge from (8, 3) to (8.3, 3.3)
\edge from (8.7, 3.7) to (9, 4)
\edge from (8, 4) to (9, 3)
\edge from (9, 1) to (10, 1)
\edge from (9, 2) to (10, 2)
\edge from (9, 3) to (10, 3)
\edge from (9, 4) to (10, 4)
\Label {$1$} at (13, 1)
\Label {$2$} at (13, 2)
\Label {$3$} at (13, 3)
\Label {$4$} at (13, 4)
\edge from (14, 1) to (15, 1)
\edge from (14, 2) to (15, 2)
\edge from (14, 3) to (15, 3)
\edge from (14, 4) to (15, 4)
\edge from (15, 1) to (16, 1)
\edge from (15, 2) to (16, 2)
\edge from (15, 3) to (15.3, 3.3)
\edge from (15.7, 3.7) to (16, 4)
\edge from (15, 4) to (16, 3)
\edge from (16, 1) to (17, 1)
\edge from (16, 2) to (16.3, 2.3)
\edge from (16.7, 2.7) to (17, 3)
\edge from (16, 3) to (17, 2)
\edge from (16, 4) to (17, 4)
\edge from (17, 1) to (17.3, 1.3)
\edge from (17.7, 1.7) to (18, 2)
\edge from (17, 2) to (18, 1)
\edge from (17, 3) to (18, 3)
\edge from (17, 4) to (18, 4)
\edge from (18, 1) to (19, 1)
\edge from (18, 2) to (19, 2)
\edge from (18, 3) to (18.3, 3.3)
\edge from (18.7, 3.7) to (19, 4)
\edge from (18, 4) to (19, 3)
\edge from (19, 1) to (20, 1)
\edge from (19, 2) to (19.3, 2.3)
\edge from (19.7, 2.7) to (20, 3)
\edge from (19, 3) to (20, 2)
\edge from (19, 4) to (20, 4)
\edge from (20, 1) to (21, 1)
\edge from (20, 2) to (21, 2)
\edge from (20, 3) to (21, 4)
\edge from (20, 4) to (20.3, 3.7)
\edge from (20.7, 3.3) to (21, 3)
\edge from (21, 1) to (22, 1)
\edge from (21, 2) to (22, 2)
\edge from (21, 3) to (22, 3)
\edge from (21, 4) to (22, 4)
\endgraph
\hfill}
\caption{Ziegler's pair: 
The half-braids $\a^{+}$ (left) and $\a^{-}$ (right).}
\label{fig:zbraids}
\end{figure}

\subsection{} \label{rigidisotop}  
Finally, let us consider the natural isotopy relation between arrangements, 
modeled on the similar notion for configurations.  
  
\begin{defn}  \label{def:isotopic}
Two arrangements $\A$ and $\A'$ are called {\em rigidly isotopic} if there is 
an isotopy of $\R^{4}$ connecting $\A$ to $\A'$ through arrangements. 
\end{defn}

The rigid isotopy class of $\CC_H(\A)$ does not 
depend on $H$, and the rigid isotopy class of $\A_p(\CC)$ does not 
depend on $p$.  Therefore, we will denote them simply by $\CC(\A)$ and 
$\A(\CC)$, respectively.  Moreover, rigid isotopy classes of configurations 
are in one-to-one correspondence with rigid isotopy classes of $2$-arrangements.  
See Crapo and Penne~\cite{CP} for details and references.  

\begin{rem} \label{rem:topmove}
Given an arrangement, we can deform it by means of a rigid isotopy 
so that one of the planes has linking number $+1$ will all other planes. 
The analogous procedure for bringing one of the lines of a configuration 
on top of all others is explained in Penne~\cite{P}. 
\end{rem}


\section{Decomposable and Horizontal Arrangements} 
\label{sec:decomphoriz}

In this section we look at arrangements that can 
be obtained by a sequence of cabling operations from simpler 
arrangements, and also at arrangements whose corresponding 
configurations are ``horizontal".  We consider in more 
detail the subclass of completely decomposable arrangements, 
and obtain a normal form for those of depth~$2$. 

\subsection{}  
\label{cables}
Let us start by recalling the following notion from knot theory 
(see~\cite{BZ}).  Let $L=L_1 \cup \cdots \cup L_n$ be a link in $\SP^3$.  
The $(a,b)$-cable of $L$ about the $k^{\text{th}}$ component 
is the link $L\{a,b\}=L\cup K(a,b)$, where $K(a,b)$ is an  
$(a,b)$-torus link contained in the boundary 
of a tubular neighborhood of $L_k$.  

Now let $\A$ be a $2$-arrangement of $n$ planes in $\R^4$, with 
defining polynomial $f=f_1\cdots f_n$.  Fix an index $1\le k\le n$, 
a positive integer $r$, and a number $\epsilon=\pm 1$.  Given these 
data, we define the {\em $\epsilon r$-cable} about the $k^{\text{th}}$ 
component of $\A$ to be the arrangement $\A^{k}\{\epsilon r\}$ 
with defining polynomial 
\begin{equation*}  \label{def:cabling}  
f(f_k+g_1)\cdots (f_k+g_r),
\end{equation*}
where each $g_j$ is a linear form in $z, \bar{z}, w, \bar{w}$, whose 
coefficients are sufficiently small with respect to those of $f$, and 
such that $\sgn(\det(f_k,g_j))=\epsilon$, for $j=1,\dots, r$. 

The cabling operation is well-defined up to rigid isotopy of arrangements.  
The reverse operation is called decabling.  It is readily seen that the 
link of $\A^{k}\{\pm r\}$ is the $(r, \pm r)$-cable about 
the $k^{\text{th}}$ component of $L(\A)$.  

\begin{defn}  \label{def:compdecomp}  
A $2$-arrangement for which no decabling is possible is called 
{\it indecomposable}; otherwise, it is called {\it decomposable}.
If $\A$ is connected to the trivial arrangement $\A_1$ by a finite 
sequence of cabling moves, then $\A$ is called {\em completely decomposable}.  
\end{defn}

\begin{exmp} \label{ex:compdecomp}
The complex arrangement $\A_n$ is the $(n-1)$-cable of $\A_1$, 
and thus is completely decomposable.  Its link is the corresponding 
cable about the unknot.  The arrangement $\A^{-}$ from Example~\ref{ex:ziegler} 
is the $(-1)$-cable of $\A_3$, and thus is also completely decomposable.  
\end{exmp}

\subsection{}  
We now define $2$-arrangements in $\R^4$ corresponding to 
special collections of skew lines in $\R^3$, variously called 
join configurations \cite{V}, horizontal configurations \cite{M1}, 
or spindle configurations \cite{CP}.

\begin{defn}
A configuration is called {\em horizontal} if it is rigidly isotopic 
to a configuration whose lines are stacked one over another in distinct  
planes, all parallel to a fixed (horizontal) plane.  A $2$-arrangement 
which admits an associated horizontal configuration is called {\em horizontal}.
\end{defn}

A horizontal configuration $\CC$ of $n$ lines determines a permutation
$\tau=\tau(\CC)$ on $\{1,\dots , n\}$, as follows.  Project perpendicularly 
all lines onto a fixed horizontal plane.  Order these $n$ lines 
in decreasing order of their (necessarily distinct) slopes.  
Order the $n$ horizontal planes containing the lines in increasing order 
of their vertical heights.  For every $i\in\{1,\dots, n\}$, put $\t_{i}=k$ 
if the $i^{\text{th}}$ line is contained in the $k^{\text{th}}$ horizontal plane.  
This defines the permutation $\tau\in S_n$.  

Conversely, every permutation $\t\in S_n$ determines a horizontal 
configuration $\CC(\t)$ (see \cite{DV, M1}), and thereby a 
horizontal arrangement $\A(\t)$.  Explicitly,  $\A(\t)$ may 
be defined as follows.  

\begin{prop} \label{prop:horizontal} 
Let $\t\in S_n$. Choose real numbers 
$a_i, b_i$, $1\le i\le n$, so that
$a_1 < \dots < a_n$ and $b_{\t_{1}} < \dots < b_{\t_{n}}$.  
Then the polynomial 
\begin{equation}  \label{eq:horizpoly}
f(z,w)=\prod_{i=1}^{n}{(z-\frac{a_i+b_i}{2}w-\frac{a_i-b_i}{2}\bar{w})},
\end{equation}
defines a horizontal $2$-arrangement, whose associated 
permutation is $\t$. 
\end{prop}

For horizontal arrangements, the linking numbers have a 
particularly simple interpretation.  Namely, if $\A=\A(\t)$, 
then $l_{i,j}=\sgn(\t_{i}\t_{j})$.  

\begin{exmp}  \label{ex:zieghoriz}
In Example~\ref{ex:ziegler}, pick the vertical coordinate to be 
$y=\Im z$. Then the lines of the configurations $\CC^{\pm}$ are 
contained in horizontal planes, parallel to the plane $y=0$. 
In each case, the ordering given by the slopes is $(1,2,3,4)$.  
The ordering given by the vertical heights is 
$(1,2,3,4)$ for $\CC^{+}$, and $(2,1,3,4)$ for $\CC^{-}$.  
Thus $\A^{+}=\A(1234)$ and $\A^{-}=\A(2134)$.  The defining 
polynomials corresponding to the choices $a=(1,2,3,4)$, 
$b=(\pm 1, \pm 2, 3, 4)$ in \eqref{eq:horizpoly} are 
the polynomials $f^{\pm}$ from Example~\ref{ex:ziegler}.  
All linking numbers $l_{i,j}^{\pm}$ are equal 
to $+1$, except for $l_{1,2}^{-}=-1$.  
\end{exmp}  

The permutation $\t$ associated to a horizontal arrangement $\A=\A(\t)$ 
is not unique.  The following result of Mazurovski\u{\i}~\cite{M2} 
lists various ways in which uniqueness is known to fail.  

\begin{prop}  
\label{prop:maz} 
Two horizontal arrangements, defined by permutations $\t$ and $\t'$ in $S_n$, 
are rigidly isotopic if:  
\begin{alphenum}
\item \label{move-a} $\t'=\sigma\t\sigma'$, where 
$\sigma$ and $\sigma'$ are circular permutations of $(1,\dots , n)$;  or 
\item \label{move-b} $\t'=\t^{-1}$;  or
\item \label{move-c}
$\t'=(\t_1, \dots, \t_{i-1}, \t'_{i}, \dots, \t'_{i+s}, \t_{i+s+1}, \dots, \t_n)$, where 
$(\t_{i}, \dots, \t_{i+s})$ is a permutation of $\{m+1,\dots ,m+s+1\}$, and 
 $(\t'_{i}-m, \dots ,\t'_{i+s}-m)=(s+1, \dots, 1)(\t_{i}-m, \dots, 
\t_{i+s}-m)(s+1, \dots, 1)$.  
\end{alphenum}
\end{prop}

\begin{rem} \label{rem:horizconj}
We do not know whether any two rigidly isotopic horizontal arrangements can 
be connected by a finite sequence of moves of type (\ref{move-a}), (\ref{move-b}), 
(\ref{move-c}).  There is another set of moves, introduced by Crapo and Penne, 
which is conjectured to be complete for horizontal configurations, see~\cite{CP}, p.~80.  
At any rate, the precise enumeration of the cosets of $S_n$  modulo the 
equivalence relation generated by either set of moves seems to be a 
challenging combinatorial problem.  
\end{rem}

\begin{exmp} \label{rem:fixmove}
We can use moves of type~(\ref{move-a}) to realize the rigid isotopy 
from Remark~\ref{rem:topmove} in the case of horizontal arrangements.  
Indeed, if $\A=\A(\tau)$ for some $\t\in S_n$ with $\t_k=n$, then we 
can replace $\t$ by $\t'=\t(k+1,\dots , n,1,\dots , k)$.  
This yields a new arrangement, $\A'=\A(\t')$, for which 
$\t'_n=n$, and $l'_{i,n}=1$ for all $i<n$.  
\end{exmp}

\begin{exmp} \label{rem:blockmove}
An important example of move~(\ref{move-c}) is as follows.   
Suppose the block $B=(\t_i,\dots, \t_j)$ is obtained  by concatenation 
of two blocks, $B_1$ and $B_2$, of consecutive integers,  each block 
in either increasing order (a {\em positive block}), or in decreasing 
order (a {\em negative block}), and so that $\min B_2=\max B_1+1$.  
Then the block $B'=(\t'_i, \dots, \t'_j)$ is also a concatenation 
of two blocks  of consecutive integers, $B'_1$ and $B'_2$. Moreover, 
$B'_1$ is $B_2$  shifted down by $|B_1|$ and $B'_2$ is $B_1$ shifted up by
$|B_2|$. In essence, the move $(B_1 B_2)\to(B'_1 B'_2)$ allows us to flip-and-shift
adjacent blocks of consecutive integers, provided that $\min B_2=\max B_1+1$.
\end{exmp}

\subsection{}  \label{subsec:shift}
We come now to a special class of horizontal arrangements, 
that can be constructed inductively from $\A_1$ by a sequence 
of cabling moves.  Let $\A=\A(\t)$, where $\t\in S_n$, and $k=\t_{\ell}$.  
An $\epsilon r$-cabling move on the $k^{\text{th}}$ component of $\A$ 
yields a new horizontal arrangement, $\A(\t')$, where $\t'\in S_{n+r}$ 
is given by
\begin{equation*}  \label{eq:shift}
\t'_{i}=
\begin{cases}
\t_{i} & \text{if } i\notin\{\ell,\dots,\ell+r\} \text{ and } \t_{i}<k\\
\t_{i}+r& \text{if } i\notin\{\ell,\dots,\ell+r\} \text{ and } \t_{i}>k\\
k+\epsilon(\ell-i)+\frac{1-\epsilon}{2} r& \text{if } i\in\{\ell,\dots,\ell+r\}.
\end{cases}
\end{equation*}
In other words, an $\epsilon r$-cabling move on $k$
shifts all the numbers in $\t$ greater than $k$ by $r$ and replaces $k$ by
$(k,\dots, k+r)$ if $\epsilon=1$ or by $(k+r,\dots, k)$ if $\epsilon=-1$.

\begin{defn}  
A horizontal arrangement is called {\em completely decomposable} if 
the associated permutation can be obtained from the
identity permutation $(1)$ by a finite sequence of cabling moves.  
\end{defn} 

It is apparent from the definitions that the link of 
a completely decomposable arrangement is obtained from the unknot 
by successive $(1,\pm 1)$-cablings. 

\begin{exmp}  \label{ex:decomp} 
All $2$-arrangements of up to $5$~planes are completely decomposable, 
except for $\A(31425)$, which is indecomposable.  Among arrangements 
of $6$~planes, for example, $\A(K)=\A(341256)$ is completely decomposable, 
$\A(314256)$ is decomposable but not completely so, and $\A(241536)$ 
is indecomposable. 
\end{exmp}

\subsection{}
We now introduce a measure of the complexity of a completely 
decomposable arrangement $\A$.  Pick a permutation $\t$ 
such that $\A=\A(\t)$.  Construct a sequence of permutations 
connecting $\t$ to the identity permutation, 
$\t=\t_0 \to \t_1 \to \cdots \to \t_d=(1)$, as follows.  At each step, 
partition the current permutation into blocks of consecutive integers, 
either in increasing order (positive cablings), or in decreasing order
(negative cablings), and contract each block to a single number 
(via decabling moves), renumbering accordingly.   Let $d(\t)=d$.  

\begin{defn}  \label{def:depth}
The {\em depth} of a completely decomposable arrangement $\A$ is 
\[
\depth(\A) = \min_{\t}\, \{d(\t)\mid \A \text{ is rigidly isotopic to } \A(\t)\}.
\]
\end{defn}

\begin{exmp}  \label{ex:depth}  
The only arrangement of depth~$0$ is the trivial arrangement 
$\A_1=\A(1)$.   The arrangements of depth~$1$ are 
the complex arrangements $\A_n=\A(1 \cdots n)$, with $n>1$.   
The arrangement $\A(21435)$ is completely decomposable, via 
the sequence of permutations $(21435)\to (123)\to (1)$, 
and so has depth $2$.  
\end{exmp}

For arrangements of depth~$2$, we single out the following type.   

\begin{defn} \label{defn:normalform}
Let $\A(\t)$ be a depth~$2$, completely decomposable arrangement.  
We say that $\A(\t)$ is in {\em normal form} if $\t=(I_1, \dots, I_r, J)$, 
where $I_1, \dots, I_r$ are negative blocks, $J$ is a positive block, 
and the following conditions hold:
\begin{romenum}
\item $I_1 < \cdots < I_r < J$,   \label{condition1}
\item $2 \leq |I_1| \leq \cdots \leq |I_r|$,   \label{condition2}
\item $|I_1| \le |J|$ if $r=1$.  \label{condition3}
\end{romenum} 
\end{defn}

\begin{prop}  \label{prop:normalform}  
Every arrangement of depth~$2$ is rigidly isotopic to a unique 
arrangement in normal form.   
\end{prop}

\begin{proof}
Let $\A$ be a depth~$2$ arrangement of $n$ planes.  Up to rigid isotopy, we 
may assume that $\A=\A(\t)$, where $d(\t)=2$.  Applying the type~(\ref{move-a}) 
move of Example~\ref{rem:fixmove}, we may further assume that $n$ is fixed 
by $\t$. Then the permutation sequence of $\A=\A(\tau)$ has the following form:  
$\t\to (1,\dots, r)\to(1)$.  Applying  repeatedly the type~(\ref{move-c}) move of 
Example~\ref{rem:blockmove}, we can push all the positive blocks of $\t$ 
(including singletons) to the right, packing all of them into a single 
positive block (that will contain $n$), and also arrange the negative blocks 
in increasing order of their sizes from left to the right.  In this way, we 
arrive at the normal form $\A(I_1, \dots, I_r, J)$ for $\A$.  The uniqueness 
is guaranteed by the conditions imposed on $I_1, \dots, I_r$ and $|J|$.
\end{proof}

Thus, we may refer to {\em the} normal form of an arrangement of depth~$2$.  
As we shall see in \S\ref{sec:bottomvar}, the normal form is a complete 
homotopy type invariant for complements of such arrangements.


\section{Braids and Fundamental Groups} 
\label{sec:braids}

In this section, we associate to a $2$-arrangement of $n$ planes 
several braids on $n$ strings, and use these braids to find 
presentations for the fundamental group of the complement.  

\subsection{}  \label{braids} 
Let $B_n$ be Artin's braid group on $n$ strings, with generators 
$\s_1, \s_2, \dots , \s_{n-1}$ and relations $\s_i\s_{j}\s_i=
\s_{j}\s_i\s_{j}$ for $|i-j|=1$ and $\s_i\s_j=\s_j\s_i$ for $|i-j|>1$, 
see~\cite{Bi}.  Also, let 
$\D_n=(\s_{n-1} \cdots \s_1)(\s_{n-1} \cdots \s_2) \cdots
(\s_{n-1}\s_{n- 2})(\s_{n-1})\in B_n$ be 
``Garside's braid"---the half-twist on $n$ strings.

Consider a configuration $\CC=\{\ell_1,\dots ,\ell_n\}$  
of $n$ skew-lines in $\R^3$.   Associated to $\CC$, 
there is a braid on $n$ strings, $\a=\a(\CC)\in B_n$, see 
Mazurovski\u{\i}~\cite{M1} and Crapo and Penne~\cite{CP}.  
The procedure that takes $\CC$ to $\a$ 
is illustrated in Figures~\ref{fig:zpair} and \ref{fig:zbraids}.  
Set $\b=\a\D_n\a\D_n^{-1}$.  We call $\a$ and $\b$, the {\em half-braid}, 
respectively the {\em full-braid} of the configuration $\CC$ 
(or of the arrangement $\A=\A(\CC)$).  As is well-known, 
conjugation by $\D_n$ is the involution $\s_i \mapsto \s_{n-i}$. Thus, 
the braid $\b$ is obtained by concatenating $\a$ with another 
copy of $\a$, rotated by $180^{\circ}$, see Figure~\ref{fig:braids}.  
Clearly, $\b$ is a pure braid in $P_n$.

\begin{figure}
\[
\beginsmallgraph
\edge from (0, 2) to (1, 2)
\edge from (0, 3) to (1, 3)
\edge from (0, 4) to (1, 4)
\edge from (1, 2) to (2, 2)
\edge from (1, 3) to (1.3, 3.3)
\edge from (1.7, 3.7) to (2, 4)
\edge from (1, 4) to (2, 3)
\edge from (2, 2) to (2.3, 2.3)
\edge from (2.7, 2.7) to (3, 3)
\edge from (2, 3) to (3, 2)
\edge from (2, 4) to (3, 4)
\edge from (3, 2) to (4, 2)
\edge from (3, 3) to (4, 4)
\edge from (3, 4) to (3.3, 3.7)
\edge from (3.7, 3.3) to (4, 3)
\edge from (4, 2) to (5, 2)
\edge from (4, 3) to (5, 3)
\edge from (4, 4) to (5, 4)
\edge from (8, 2) to (9, 2)
\edge from (8, 3) to (9, 3)
\edge from (8, 4) to (9, 4)
\edge from (9, 2) to (10, 2)
\edge from (9, 3) to (9.3, 3.3)
\edge from (9.7, 3.7) to (10, 4)
\edge from (9, 4) to (10, 3)
\edge from (10, 2) to (10.3, 2.3)
\edge from (10.7, 2.7) to (11, 3)
\edge from (10, 3) to (11, 2)
\edge from (10, 4) to (11, 4)
\edge from (11, 2) to (12, 2)
\edge from (11, 3) to (12, 4)
\edge from (11, 4) to (11.3, 3.7)
\edge from (11.7, 3.3) to (12, 3)
\edge from (12, 2) to (12.3, 2.3)
\edge from (12.7, 2.7) to (13, 3)
\edge from (12, 3) to (13, 2)
\edge from (12, 4) to (13, 4)
\edge from (13, 2) to (14, 2)
\edge from (13, 3) to (13.3, 3.3)
\edge from (13.7, 3.7) to (14, 4)
\edge from (13, 4) to (14, 3)
\edge from (14, 2) to (15, 3)
\edge from (14, 3) to (14.3, 2.7)
\edge from (14.7, 2.3) to (15, 2)
\edge from (14, 4) to (15, 4)
\edge from (15, 2) to (16, 2)
\edge from (15, 3) to (16, 3)
\edge from (15, 4) to (16, 4)
\endgraph
\]
\caption{The braids $\a$ and $\b$ associated to $\A(213)$.}
\label{fig:braids}
\end{figure}

The following result of Mazurovski\u{\i}~\cite{M1} 
and Crapo and Penne~\cite{CP} establishes the direct connection 
between the link and the braid of an arrangement.  First recall 
the classical theorem of Alexander, according to which every link 
in $\SP^3$ is isotopic to the closure of a braid (see~\cite{Bi}).  

\begin{prop} \label{prop:closebraid} 
Let $\A$ be a $2$-arrangement in $\R^4$ and $L=L(\A)$ its link.  
Let $\CC=\CC(\A)$ be the associated configuration of skew lines in $\R^3$ 
and $\b=\b(\CC)$ its full-braid.  Then $L$ is isotopic to the 
closure of $\b$.
\end{prop}

Let $X$ be the complement of the arrangement $\A$, and 
$G=\pi_1(X)$ its fundamental group.  Recall that $X$ is  
homotopy equivalent to the complement $Y$ of the link $L$.  
Since $L$ is the closure of $\b$, the group $G$ has 
Artin presentation 
\begin{equation} \label{eq:artinpres}
G=\langle x_1,\dots , x_n \mid  \b(x_i) = x_i, 
\: i=1,\dots , n\rangle, 
\end{equation}
see~\cite{Bi, BZ}. 

\subsection{}
As mentioned in Remark~\ref{rem:topmove}, we can bring one of the 
lines of $\CC$, say $\ell_n$, on top of all the other ones.  
Discarding $\ell_n$, we get a configuration $\check{\CC}$ 
of $n-1$ skew lines, so that $\CC=\CC'\cup\{\ell_n\}$.  
It follows that $L=\check{L}\cup L_n$, where $\check{L}$ 
is the closure of $\check{\b}=\b(\check{\CC})\in P_{n-1}$.  
Furthermore, it is readily seen that the half-braid of 
$\check{\CC}$ is given by 
$\iota(\check{\a})=\s_1^{-1}\cdots \s_{n-1}^{-1}\a$, 
where $\iota:B_{n-1}\hookrightarrow B_n$ is the standard 
inclusion $\iota(\s_i)=\s_{i+1}$. 
We call $\check{\a}$ and $\check{\b}$, the {\em reduced half-braid}, 
respectively the {\em reduced full-braid} of the arrangement $\A=\A(\CC)$.
   
It is now apparent that the complement of $L$ in $\SP^3$ 
is homotopy equivalent to the complement of $\check{L}$ 
in the solid torus $\SP^1\times D^2 = \SP^3 \setminus (L_n\times D^2)$.  
These geometric considerations lead to the following:

\begin{prop} \label{prop:bundle}
The complement $X(\A)$ of a $2$-arrangement of $n$ planes in $\R^4$
is homotopy equivalent to the total space of a bundle over the circle, 
with fiber $D^2\setminus \{n-1 \text{ points}\}$, and monodromy 
the braid automorphism $\check{\b}$.
\end{prop}

Thus, $X$ is a $K(G,1)$, with fundamental group a semidirect product 
of free groups, $G = \F_{n-1} \rtimes_{\check{\b}} \F_1$.
The Artin representation of $\check{\b}$ provides a presentation for $G$, 
corresponding to this split extension:
\begin{equation*} \label{eq:bpres}
G=\langle x_1,\dots , x_n \mid x_n^{-1} x_i x_n = \check{\b}(x_i), 
\: i=1,\dots , n-1\rangle. 
\end{equation*}

\begin{exmp}\label{ex:complex} 
For the complex arrangement $\A_n$, the half-braid is the 
half-twist $\a=\D_n$, and the full-braid is the full-twist 
$\beta=\D_n^2$.  Since $\check{\b}=\D_{n-1}^2$ 
acts on $\F_{n-1}$ by conjugation by $x_1\cdots x_{n-1}$, 
the group $G$ is isomorphic to $\F_{n-1}\times \F_1$, where 
$\F_1=\langle x_1\cdots x_n\rangle$.
\end{exmp}

For a non-complex $2$-arrangement, the group $G$ is in general 
{\em not} isomorphic to a direct product, as we shall later see.  
Nevertheless, we may still use the underlying idea of  
Example~\ref{ex:complex}, and simplify the presentation of $G$, 
by cutting off a full twist from $\check{\b}$.

\begin{prop}  \label{prop:xi2}
Let $\A$ be a $2$-arrangement of $n$ planes, with reduced 
half-braid $\Check{\a}$.  Set $\xi =\D_{n-1}{\check{\a}}^{-1}$.  
Then, the fundamental group $G$ of $\A$ is 
isomorphic to $\F_{n-1} \rtimes_{\xi^2} \F_1$, and has presentation
\begin{equation}  \label{eq:xipres}
G=\langle x_1,\dots , x_n \mid x_n x_i x_n^{-1} = \xi^2(x_i), 
\: i=1,\dots , n-1\rangle. 
\end{equation}
\end{prop}
\begin{proof}  
Recall that $G_1\rtimes_{\phi} G_2\cong G_1\rtimes_{\phi'} G_2$ 
if $\phi'=\gamma \phi^{\pm 1}$, where $\gamma\in\Inn(G_1)$.  
Thus, it suffices to show that $\check{b}$ differs from $\xi^{-2}$ by an 
inner automorphism of $\F_{n-1}$.  This follows from the fact that 
$\D_{n-1}^2\in \Center(B_{n-1})\cap \Inn(\F_{n-1})$:  
\begin{equation*}
\hspace{0.6in}
\check{b}= \check{a} \D_{n-1} \check{a} \D_{n-1} 
   = \xi^{-1} \D_{n-1}^2 \xi^{-1} 
   = \D_{n-1}^2 \xi^{-2}. 
\hspace{0.6in}\qed
\end{equation*}
\renewcommand{\qed}{}
\end{proof}

\begin{rem} \label{rem:conjugate}
Recall also that $G_1\rtimes_{\phi} G_2\cong G_1\rtimes_{\phi'} G_2$ 
if $\phi'=\psi^{} \phi \psi^{-1}$.  We can use this observation 
to further simplify the above presentation, 
by conjugating $\xi\in P_{n-1}$ by a suitable automorphism of $\F_{n-1}$.  
In practice, this will be achieved by either changing the basis 
of $\F_{n-1}$, or by conjugating $\xi$ by a suitable braid 
$\delta\in B_{n-1}$.  
\end{rem}

\subsection{} 
We now identify the braids $\a\in B_{n+1}$ and $\xi\in P_{n}$ 
associated to a horizontal arrangement of $n+1$ planes in terms of the 
generators $\s_i$ of the braid group and of the generators 
$A_{i,j}=\s_{j-1} \cdots \s_{i+1}\s_{i}^{2}\s_{i+1}^{-1} 
\cdots \s_{j-1}^{-1}$ of the pure braid group, respectively.    

\begin{prop}  \label{prop:combed}
Let $\A(\t)$ be a horizontal $2$-arrangement of $n+1$ planes.  Then:
\begin{alphenum} 
\item \label{al}
The half-braid $\a$ has the form
\[
\a=(\s_{n}^{l_{n,n+1}} \cdots \s_1^{l_{1,n+1}})(\s_{n}^{l_{n-1,n}} 
\cdots \s_2^{l_{1,n}}) \cdots
(\s_{n}^{l_{2,3}}\s_{n-1}^{l_{1,3}})(\s_{n}^{l_{1,2}}), 
\]
where $l_{i,j}$ is the sign of the permutation $(\t_{i}\, \t_{j})$. 
\item \label{xi}
The pure braid $\xi$ can be combed as 
$\xi=\xi_2 \cdots \xi_{n}$, where 
\begin{equation*}  \label{eq:combed}
\xi_j=\prod_{i=1}^{j-1} A_{i,j}^{e_{i,j}} \quad \text{ and } \quad
e_{i,j}=
\begin{cases}
1 &\text{ if } \t_{i} > \t_{j},\\
0 &\text{ otherwise.}
\end{cases}
\end{equation*}
\end{alphenum}
\end{prop}

\begin{proof}  
Part \eqref{al} follows from the definitions of $\CC=\CC(\t)$ and 
$\a=\a(\CC)$, and the fact that the linking numbers of $L(\A(\t))$ 
are given by $l_{i,j}=\sgn(\t_{i}\t_{j})$.   

For part \eqref{xi}, it is enough to show that, for all 
$k$ with $1 \leq k \leq n-1$,  
\begin{equation} \label{eq:iden1}
\begin{split} 
&\D_n (\s_{n-1}^{-l_{1,2}}) (\s_{n-2}^{-l_{1,3}}\s_{n-1}^{-l_{2,3}}) \cdots
(\s_{n-k}^{-l_{1,k+1}} \cdots \s_{n-1}^{-l_{k,k+1}}) = \\
&\qquad (A_{1,2}^{e_{1,2}}) \cdots (A_{1,k+1}^{e_{1,k+1}} \cdots 
A_{k,k+1}^{e_{k,k+1}})
\ (\s_{k+1} \cdots \s_{n-1})\cdots (\s_2 \cdots \s_{n-k}) \D_{n-k} ,  
\end{split}
\end{equation}
Indeed, the identity~\eqref{eq:iden1} for $k=n-1$ yields 
the desired combed form of $\xi$.

The proof of \eqref{eq:iden1} is by induction on $k$, using the 
braid relations.  The step $k=1$ is as follows: 
\begin{align*}  
\D_n \s_{n-1}^{\pm 1} 
&=(\s_1 \cdots \s_{n-1}) \cdots (\s_1\s_2)(\s_1) \s_{n-1}^{\pm 1} \\
&=(\s_1)(\s_2\s_1) \cdots (\s_{n-2}\s_{n-3})(\s_{n-1}\s_{n-2}) \s_{n-1}^{\pm 1}
(\s_1 \cdots \s_{n-3}) \cdots (\s_1\s_2)(\s_1) \\
&=(\s_1)(\s_2\s_1) \cdots (\s_{n-2}\s_{n-3}) \s_{n-2}^{\pm 1} (\s_{n-1}\s_{n-2})
(\s_1 \cdots \s_{n-3}) \cdots (\s_1\s_2)(\s_1)  \\
&\quad\cdots \\
&=(\s_1) \s_1^{\pm 1} (\s_2\s_1) \cdots (\s_{n-2}\s_{n-3})(\s_{n-1}\s_{n-2})
(\s_1 \cdots \s_{n-3})  \cdots  (\s_1\s_2)(\s_1) \\
&=A_{1,2}^{e} (\s_2 \cdots \s_{n-1}) \D_{n-1}, 
\end{align*}
where $e=\frac{1\pm 1}{2}$. The induction step is similar but tedious, 
and will be omitted. 
\end{proof}

\begin{figure}
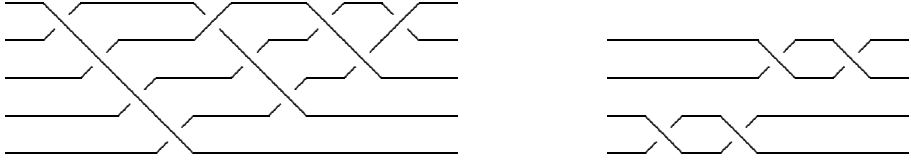

\[
\beginsmallgraph
\edge from (2, 1) to (3, 1)
\edge from (2, 2) to (3, 2)
\edge from (2, 3) to (3, 3)
\edge from (2, 4) to (3, 4)
\edge from (2, 5) to (3, 5)
\edge from (3, 1) to (4, 1)
\edge from (3, 2) to (4, 2)
\edge from (3, 3) to (4, 3)
\edge from (3, 4) to (3.3, 4.3)
\edge from (3.7, 4.7) to (4, 5)
\edge from (3, 5) to (4, 4)
\edge from (4, 1) to (5, 1)
\edge from (4, 2) to (5, 2)
\edge from (4, 3) to (4.3, 3.3)
\edge from (4.7, 3.7) to (5, 4)
\edge from (4, 4) to (5, 3)
\edge from (4, 5) to (5, 5)
\edge from (5, 1) to (6, 1)
\edge from (5, 2) to (5.3, 2.3)
\edge from (5.7, 2.7) to (6, 3)
\edge from (5, 3) to (6, 2)
\edge from (5, 4) to (6, 4)
\edge from (5, 5) to (6, 5)
\edge from (6, 1) to (6.3, 1.3)
\edge from (6.7, 1.7) to (7, 2)
\edge from (6, 2) to (7, 1)
\edge from (6, 3) to (7, 3)
\edge from (6, 4) to (7, 4)
\edge from (6, 5) to (7, 5)
\edge from (7, 1) to (8, 1)
\edge from (7, 2) to (8, 2)
\edge from (7, 3) to (8, 3)
\edge from (7, 4) to (8, 5)
\edge from (7, 5) to (7.3, 4.7)
\edge from (7.7, 4.3) to (8, 4)
\edge from (8, 1) to (9, 1)
\edge from (8, 2) to (9, 2)
\edge from (8, 3) to (8.3, 3.3)
\edge from (8.7, 3.7) to (9, 4)
\edge from (8, 4) to (9, 3)
\edge from (8, 5) to (9, 5)
\edge from (9, 1) to (10, 1)
\edge from (9, 2) to (9.3, 2.3)
\edge from (9.7, 2.7) to (10, 3)
\edge from (9, 3) to (10, 2)
\edge from (9, 4) to (10, 4)
\edge from (9, 5) to (10, 5)
\edge from (10, 1) to (11, 1)
\edge from (10, 2) to (11, 2)
\edge from (10, 3) to (11, 3)
\edge from (10, 4) to (10.3, 4.3)
\edge from (10.7, 4.7) to (11, 5)
\edge from (10, 5) to (11, 4)
\edge from (11, 1) to (12, 1)
\edge from (11, 2) to (12, 2)
\edge from (11, 3) to (11.3, 3.3)
\edge from (11.7, 3.7) to (12, 4)
\edge from (11, 4) to (12, 3)
\edge from (11, 5) to (12, 5)
\edge from (12, 1) to (13, 1)
\edge from (12, 2) to (13, 2)
\edge from (12, 3) to (13, 3)
\edge from (12, 4) to (13, 5)
\edge from (12, 5) to (12.3, 4.7)
\edge from (12.7, 4.3) to (13, 4)
\edge from (13, 1) to (14, 1)
\edge from (13, 2) to (14, 2)
\edge from (13, 3) to (14, 3)
\edge from (13, 4) to (14, 4)
\edge from (13, 5) to (14, 5)
\edge from (18, 1) to (19, 1)
\edge from (18, 2) to (19, 2)
\edge from (18, 3) to (19, 3)
\edge from (18, 4) to (19, 4)
\edge from (19, 1) to (19.3, 1.3)
\edge from (19.7, 1.7) to (20, 2)
\edge from (19, 2) to (20, 1)
\edge from (19, 3) to (20, 3)
\edge from (19, 4) to (20, 4)
\edge from (20, 1) to (21, 1)
\edge from (20, 2) to (21, 2)
\edge from (20, 3) to (21, 3)
\edge from (20, 4) to (21, 4)
\edge from (21, 1) to (21.3, 1.3)
\edge from (21.7, 1.7) to (22, 2)
\edge from (21, 2) to (22, 1)
\edge from (21, 3) to (22, 3)
\edge from (21, 4) to (22, 4)
\edge from (22, 1) to (23, 1)
\edge from (22, 2) to (23, 2)
\edge from (22, 3) to (22.3, 3.3)
\edge from (22.7, 3.7) to (23, 4)
\edge from (22, 4) to (23, 3)
\edge from (23, 1) to (24, 1)
\edge from (23, 2) to (24, 2)
\edge from (23, 3) to (24, 3)
\edge from (23, 4) to (24, 4)
\edge from (24, 1) to (25, 1)
\edge from (24, 2) to (25, 2)
\edge from (24, 3) to (24.3, 3.3)
\edge from (24.7, 3.7) to (25, 4)
\edge from (24, 4) to (25, 3)
\edge from (25, 1) to (26, 1)
\edge from (25, 2) to (26, 2)
\edge from (25, 3) to (26, 3)
\edge from (25, 4) to (26, 4)
\endgraph
\]
\caption{The braids $\a$ and $\xi$ associated to $\A(21435)$.} 
\label{fig:alphaxi}
\end{figure}

\begin{exmp} 
\label{ex:horizbraids}  
The complex arrangement $\A_n$ is horizontal, with corresponding 
permutation the identity $\tau=(1\cdots n)$.  Thus 
$\a=\Delta_n$ and $\xi=1$.  The arrangement $\A=\A^{-}$ 
from Example~\ref{ex:ziegler} is also horizontal, 
with permutation $\tau=(2134)$.  Thus $\a=\s_3\s_2\s_1\s_3\s_2\s_3^{-1}$ 
and $\xi=A_{1,2}$.  The braids 
$\a=\s_4\s_3\s_2\s_1\s_4^{-1}\s_3\s_2\s_4\s_3\s_4^{-1}$ 
and $\xi=A_{1,2}A_{3,4}$ associated to the horizontal arrangement 
$\A(21435)$ are illustrated in Figure~\ref{fig:alphaxi}.  
\end{exmp}


\section{Determinantal Ideals and Characteristic Varieties} 
\label{sec:detideals}

We start this section with a review of the determinantal ideals 
of the Alexander module of a space, following Hillman~\cite{Hi} and 
Turaev~\cite{T1, T2}.  From the varieties defined by 
these ideals, we extract numerical homotopy-type invariants, 
that will be used for the rest of this paper.  

\subsection{}
Let $X$ be a connected, finite CW-complex, with basepoint $*$, and 
fundamental group $\pi_1(X,*)$.  Let $p:\widetilde{X}\to X$ be the 
universal abelian cover, corresponding to the abelianization 
homomorphism $\ab:\pi_1(X,*)\to H_1(X;\Z)$. The relative 
homology group $A(X)=H_1(\widetilde{X},p^{-1}(*); \Z)$ has 
the structure of a (left) module over the group ring $\Z[H_1(X; \Z)]$, 
and is known as the {\em Alexander module} of $X$. 

Now assume that $H_1(X,\Z)$ is isomorphic to 
$\Z^n$, the free abelian group on $t_1,\dots ,t_n$.  
A choice of isomorphism, $\psi: H_1(X) \xrightarrow{\simeq} \Z^n$, 
identifies $\Z[H_1(X)]$ with $\L=\Z[t_1^{\pm 1}, \dots , t_n^{\pm 1}]$, 
the ring of Laurent polynomials in $n$ variables, and defines a 
$\Lambda$-module structure on the Alexander module of $X$, 
which we will denote by $A(X,\psi)$.  From a presentation of 
the fundamental group, 
$\pi_1(X)=\langle x_1,\dots , x_q \mid r_1,\dots , r_s\rangle$, 
one gets a presentation for the Alexander module, 
\[
\L^s \xrightarrow{M} \L^q\to A(X,\psi)\to 0,
\]
where $M=\begin{pmatrix} \partial r_i/ \partial x_j \end{pmatrix}^{\ab}$ 
is the abelianized Jacobian matrix of Fox derivatives. 

Define the {\em $k^{\text{th}}$ determinantal ideal} of $A_{\psi}(X)$ 
to be the ideal $E_k(X,\psi)$ generated by the codimension~$k$ minors 
of the Alexander matrix $M$.  Clearly, $E_{k}(X,\psi)\subseteq E_{\ell}(X,\psi)$ 
if $k\le \ell$.  The determinantal ideals depend only on the homotopy 
type of $X$ (in fact, only on its fundamental group), and on the 
identification $\psi: H_1(X)\to\Z^n$. 

If $\pi_1(X)$ has positive deficiency (i.e., admits a presentation 
with more generators than relations), then $E_1(X,\psi)$ is of the form 
$I \cdot (\D_{X,\psi})$, where $I$ is the augmentation ideal of $\L$, 
and $\D_{X,\psi}\in \L$ is the (multi-variable) {\em Alexander polynomial} 
of $X$, see~\cite{EN}. 

\subsection{}
We now associate to $X$ subvarieties $V_k(X,\psi)$ of the algebraic 
torus $(\C^{*})^n$, defined by the determinantal ideals $E_{k}(X,\psi)$, 
following \cite{DF, L1}.  
The coordinate ring of $(\C^{*})^n$ is $\L_{\C}=\L\otimes \C$, the ring 
of Laurent polynomials with complex coefficients.  Then, for each $k\ge 0$, 
we set 
\[
V_k(X,\psi)=\{  (t_1,\dots ,t_n) \in (\C^*)^n \mid g(t_1,\dots, t_n)=0, 
\text{ for all } g\in \sqrt{E_{k}(X,\psi)\otimes \C}  \},
\] 
where $\sqrt{\mathfrak{a}}$ denotes the radical of an ideal $\mathfrak{a}$. 
Clearly, $V_{k}(X,\psi)\supseteq V_{\ell}(X,\psi)$ if $k\le \ell$.

\begin{defn}    \label{def:monomisom}
Two algebraic subvarieties $V$ and $V'$ of $(\C^{*})^n$ 
are said to have the same {\it monomial isomorphism type} if there exists 
an automorphism $\f_A:(\C^{*})^n\to (\C^{*})^n$ of the form 
\[
\f_A(t_i)=t_1^{a_{i_1}} \cdots t_n^{a_{i_n}}, \quad 1 \leq i \leq n, 
\]
for some matrix $A=(a_{i,j})\in \GL_n(\Z)$, which maps $V$ into $V'$.
\end{defn}

\begin{prop}   \label{prop:towerinvariant}
The monomial isomorphism type of the subvariety $V_k(X,\psi)$ of the 
algebraic torus $(\C^{*})^n$ depends only on the isomorphism type of 
$\pi_1(X)$, and not on the identification $\psi:H_1(X)\to \Z^n$.  We call 
$V_k(X)=V_k(X,\psi)$ the {\em $k^{\text{th}}$ characteristic variety} of $X$.  
\end{prop}
\begin{proof}
Let $X$ and $Y$ be connected, finite CW-complexes, and let 
$h:\pi_1(X) \to \pi_1(Y)$ be an isomorphism.  Let 
$h_*:H_1(X) \to H_1(Y)$ be the abelianization of $h$, 
and set $\bar{h}=\psi_{Y}^{} h_* \psi_{X}^{-1}:\Z^n\to \Z^n$.  
The extension of $\bar{h}$ to $\L_{\C}=\C\Z^n$ restricts 
to an isomorphism $E_k(X,\psi_X)\otimes\C \to E_k(Y,\psi_Y)\otimes\C$, 
for each $k\ge 0$.  Now let $\f$ the (monomial) automorphism 
of $(\C^{*})^n$ induced by $\bar{h}$.  Clearly, $\f$ restricts 
to an isomorphism $V_k(X,\psi_X) \to V_k(Y,\psi_Y)$.
\end{proof}

In other words, for each $k\ge 0$, the monomial isomorphism type 
of $V_k(X)$ is an isomorphism type of $\pi_1(X)$, and thus, a 
homotopy-type invariant of $X$.  Furthermore, if $\pi_1(X)$ has 
positive deficiency, the Alexander polynomial $\Delta_{X}=\Delta_{X,\psi}$ 
is well-defined up to a monomial change of basis in $(\C^{*})^n$, 
and up to multiplication by a unit 
$c t_1^{i_1}\cdots t_n^{i_n}\in \Lambda_{\C}$.  Note that 
$V_1(X)=\mathbf{1}\cup \{\Delta_{X}=0\}$, where $\mathbf{1}=(1,\dots ,1)$ 
is the origin of the complex torus $(\C^*)^n$.  

\subsection{} 
By themselves, the characteristic varieties are not very practical homotopy-type 
invariants.  We extract from them several numerical invariants that 
are powerful enough for our purposes.  For each integer $p\ge 2$, 
let
\[\Omega_p^n=\{(\omega_1,\dots ,\omega_n)\in (\C^*)^n \mid 
\omega_i\: \text{ is a $p^{\text{th}}$ root of unity}\}\] 
be the set of $p$-torsion points of $(\C^*)^n$. 

\begin{thm} \label{thm:invariant} 
The following are isomorphism type invariants of $\pi_1(X)$:
\begin{alphenum}
\item \label{inv1}
The list $\Sigma_k(X)$ of codimensions of irreducible 
components of $V_k(X)$; 

\item \label{inv2}
The list $\Sigma_{\mathbf{1},k}(X)$ of codimensions of irreducible 
components of $V_k(X)$ passing through $\mathbf{1}$;

\item \label{inv3}
The number $\Tors_{p,k}(X)=\left|{\Omega_p^n\cap V_k(X)}\right|$ 
of $p$-torsion points of $V_k(X)$.
 
\end{alphenum}
\end{thm}

\begin{proof} 
By Proposition~\ref{prop:towerinvariant}, an isomorphism of fundamental 
groups determines a monomial isomorphism of the corresponding 
characteristic varieties.  Part~\eqref{inv1} follows from the fact 
that an isomorphism of algebraic varieties sends irreducible 
components to irreducible components of the same codimension.   
Part~\eqref{inv2} follows from Part~\eqref{inv1}, and the fact that 
a monomial isomorphism fixes $\mathbf{1}$.  Part~\eqref{inv3} follows 
from the fact that a monomial isomorphism preserves the set of 
$p$-torsion points.  
\end{proof}

\subsection{}  \label{alexmod}
Now let $X=X(\A)$ be the complement of a $2$-arrangement of 
$n$ planes in $\R^4$.  Recall that $X$ has the homotopy type 
of a $2$-complex (modeled on the Artin presentation of its 
fundamental group $G$), and that $H_1(X)=\Z^n$.  Thus, we 
can define the {\em $k^{\text{th}}$ characteristic variety of $\A$} 
to be $V_k(\A)=V_k(X)$.  As we shall see, the descending tower of 
characteristic varieties has the form $(\C^{*})^n=V_0 \supseteq V_1 \supseteq 
\dots \supseteq V_{n-2} \supseteq V_{n-1}\supseteq V_{n}=\emptyset$, 
with $V_1$ being a hypersurface in $(\C^{*})^n$, if $n\ge 3$, 
and $V_{n-1}$ consisting of the single point $\mathbf{1}$, if $n\ge 2$.  
We will focus on the nontrivial ends of the tower, namely 
$V_{1}$ and $V_{n-2}$, which we shall call the {\em top}, 
respectively the {\em bottom} characteristic variety of $\A$.  

In order to find explicit equations for the characteristic varieties, 
we need to choose a particular presentation for $G=\pi_1(X)$.  
Unless otherwise specified, we shall use 
the presentation \eqref{eq:xipres} associated to the semidirect 
product structure $G=\F_{n-1}\rtimes_{\xi^2} \F_1$ from 
Proposition~\ref{prop:xi2}.  This presentation yields an 
identification $\psi_{\xi}:H_1(X)\to \Z^n$.  
Let $A=A(X,\psi_{\xi})$ be the corresponding Alexander module.  
A presentation matrix for $A$ is the $(n-1)\times n$ 
(Alexander) matrix 
\begin{equation*}  
\label{alexmat}
M=\begin{pmatrix}
t_n\cdot \id-\T(\xi^2) & d_1
\end{pmatrix},
\end{equation*}
where $d_1=\left( 1-t_1 \ \cdots \ 1-t_{n-1}\right)^{\top}$ and 
$\T: P_{n-1} \to \GL_{n-1}(\L)$ is the Gassner representation 
of the pure braid group, see Birman~\cite{Bi}.  

The $k\times k$ minors of $M$ generate the determinantal ideal $E_k$, 
whose radical, $\sqrt{E_k}$, defines the $k^{\text{th}}$ 
characteristic variety $V_k=V_k(\A)$.  Note that $E_{n-1}=I$ and 
$E_{n}=\Lambda$, and so $V_{n-1}=\mathbf{1}$ and $V_{n}=\emptyset$. 

Now recall that a link group has deficiency~$1$, see e.g.~\cite{BZ}.  
Thus we may define the {\em Alexander polynomial of $\A$} to be 
$\Delta_{\A}=\Delta_{X,\psi_{\xi}}$. For $n=1$, we have 
$\Delta_{\A}=1$.  For $n>1$, we have 
\[
\Delta_{\A}(t_1,\dots ,t_n) = \frac{1}{t_n-1} 
\det \left( t_n\cdot \id-\T(\xi^2) \right),
\]  
see Penne~\cite{P}.  Thus $\Delta_{\A}=1$ for $n=2$.  
For $n\ge 3$, the triviality of 
the Gassner representation evaluated at $\mathbf{1}$ 
implies that $\mathbf{1}\in V_1(\A)$ and so 
$V_1(\A)=\{\Delta_{\A}(t_1,\dots ,t_n)=0\}$.

\begin{rem}  \label{rem:aplink}
For certain purposes, it is more natural to start from the 
Artin presentation \eqref{eq:artinpres} associated to the semidirect 
product structure $G=\F_{n-1}\rtimes_{\b} \F_1$. The resulting 
presentation, $A(X,\psi_{\b})$, for the Alexander module coincides 
with the usual presentation of the Alexander module of the  
link $L(\A)$.  We will denote the associated Alexander polynomial 
by $\Delta_{L(\A)}=\Delta_{X,\psi_{\b}}$.  
\end{rem}

\begin{exmp}

\label{ex:Cxcharvar} 
Let $\A_n$ be the arrangement of $n\ge 3$ complex lines through the 
origin of $\C^2$.  Recall that $\xi=1$ and $\beta=\D_n^2$ in this case.  
It is readily seen that $E_k=I \cdot (t_{n}-1)^{n-k-1}$.  
Thus $V_1=\dots =V_{n-2}=\{t_n-1=0\}$, and 
$\Delta_{\A_n}=(t_n-1)^{n-2}$, whereas 
$\Delta_{L(\A_n)}=(t_1\dots t_n-1)^{n-2}$.    
\end{exmp}

\begin{exmp} 
\label{ex:Zieglercharvar}
Let $\A$ be the arrangement $\A^{-}=\A(2134)$. Recall that 
$\xi=A_{1,2}$. The Artin representation of $\xi:\F_3\to \F_3$ is given 
by $\xi(x_1)=x_1x_2x_1x_2^{-1}x_1^{-1}$, $\xi(x_2)=x_1x_2x_1^{-1}$, 
$\xi(x_3)=x_3$.  Consider the new basis $y_1=x_1$, $y_2=x_1x_2$, $y_3=x_3$ 
for $\F_3$.  In this basis, $\xi(y_1)=y_2y_1y_2^{-1}$, $\xi(y_2)=y_2$, 
$\xi(y_3)=y_3$, and so the Alexander matrix is: 
\[
M= \begin{pmatrix}
t_4-t_2^2 & (t_2+1)(t_1-1) & 0 & 1-t_1 \\ 
0 & t_4-1 & 0 & 1-t_2 \\ 
0 & 0 & t_4-1 & 
1-t_3
\end{pmatrix}.
\]
The determinantal ideals are $E_1=I \cdot (t_4-1)(t_4-t_2^2)$, 
and $E_2=I \cdot (t_4-1,t_2^2-1) + \L \cdot(t_2+1)(t_1-1)(t_3-1)$. 
The characteristic varieties are
\begin{align*}  \label{eq:avz}
\begin{split}
V_1&=\{t_4-1=0\} \cup \{t_4-t_2^2=0\}, \\ 
V_2&=\{t_4-1=t_2+1=0\} \cup \{t_4-1=t_2-1=t_1-1=0\} \\ 
& \qquad\cup \{t_4-1=t_2-1=t_3-1=0\}.
\end{split}
\end{align*}

Now let $\A$ be the arrangement $\A^{+}=\A(1234)$. 
We know from Example~\ref{ex:Cxcharvar} that its characteristic 
varieties are $V_1=V_2=\{t_4-1=0\}$.  We plainly see that 
the characteristic varieties of $\A^{+}$ and $\A^{-}$ have a  
different number of components, and so are not isomorphic.  
Thus, $X^{+}\not\simeq X^{-}$.  In fact, $H^*(X^{+})\ncong H^*(X^{-})$, 
as was shown by Ziegler~\cite{Z}.  
\end{exmp}


\section{Bottom Characteristic Varieties} 
\label{sec:bottomvar}

In this section we study the bottom characteristic varieties $V_{n-2}(\A)$ 
of arrangements of $n$ planes that are obtained from the trivial 
arrangement by a sequence of cabling operations.  We obtain a 
complete characterization of these varieties when the sequence 
has length~$2$.  

\subsection{}
Let $\A$ be a depth~$2$ arrangement, and $\A(I_1,\dots ,I_r,J)$ its 
normal form, as introduced in Definition~\ref{defn:normalform}.  
By Proposition~\ref{prop:combed}, 
the pure braid $\xi$ can be taken to be $\xi=A_{I_1} \cdots A_{I_r}$, 
where $A_I$ is the full twist on the strings $I$, with the convention 
that whenever a negative block appears as a subindex of a braid generator, 
it will be understood as a set of integers in increasing order.   

Let us first consider a full twist $A_I$ corresponding to the 
indexing set $I=\{j, j+1, \dots, k\}$.  
The Artin representation of $A_I$ is given by $A_I(x_i)=x_Ix_ix_I^{-1}, i \in I$ 
and $A_I(x_i)=x_i$ otherwise, where $x_I=x_j \cdots x_k$. 
The change of basis $y_k=x_I$ and $y_i=x_i$ for $i \neq k$  
simplifies $A_I$ to $A_I(y_i)=y_ky_iy_k^{-1}, i \in I$ and $A_I(y_i)=y_i$ 
otherwise.  The Gassner representation of $A_I^2$ in this new basis 
is given by:

\begin{equation}  \label{eq:gasstwist}
\T(A_I^2)= 
\begin{pmatrix}
t_k^2  & \cdots & 0       & (t_k+1)(1-t_j)     \\ 
\vdots & \ddots &  \vdots & \vdots             \\
0      & \cdots & t_k^2   & (t_k+1)(1-t_{k-1}) \\ 
0      & \cdots & 0       & 1
\end{pmatrix}
\end{equation}

Now notice that the factors $A_{I_1}, \dots, A_{I_r}$ of $\xi$ braid on 
mutually disjoint groups of strings.  Therefore, we can change 
the basis in the free group for each block separately, as above.  
Hence, the Gassner representation $\T(\xi^2)$ is a block-diagonal matrix, 
with blocks as in~\eqref{eq:gasstwist}.  The Alexander matrix is:

\begin{equation}  \label{eq:alexdepth2}
M= 
\begin{pmatrix}
t_n-\T(A_{I_1}^2)  & \cdots & 0                  & 0               & d_1(I_1)  \\ 
\vdots             & \ddots &  \vdots            & \vdots          & \vdots     \\
0                  & \cdots & t_n-\T(A_{I_r}^2)  & 0               & d_1(I_r)  \\ 
0                  & \cdots & 0                  & t_n-\id_{|J|-1} & d_1(J)
\end{pmatrix}
\end{equation}
where $d_1(I)$ is the column vector whose entries are $t_i-1$, 
for $i\in I$.  

The radical $\sqrt{E_{n-2}}$ of the ideal of $2 \times 2$-minors of $M$ 
is generated by  
\[
t_n-1; \quad t_{k_p}^2 -1; \quad (t_{k_p}+1)(t_{i_p}-1)(t_{j_p}-1),
\]
where $k_p=\max I_p$,  for all $i_p \in I_p \setminus \{k_p\}$, 
$j_p \not\in I_p$, and $1\le p \le r$.  

The ideal $\sqrt{E_{n-2}}$ defines the characteristic variety $V_{n-2}$.  
In order to describe this subvariety of $(\C^*)^n$, we need some notation.   
Given a subset $I$ of $[1,n]=\{1,\dots,n\}$, let $\Bar{I} = [1,n]\setminus I$, 
and $\Check{I} = \Bar{I} \cup \{\max I\}$.  Also, let $T(I)$ be the 
subtorus of $(\C^*)^n$ given by $T(I)=\{t_i-1=0 \mid i\in I\}$, 
and $\Bar{T}(I)$ the translated subtorus given by 
$\Bar{T}(I)=\{t_i+1=0 \mid i\in I\}$.  

\begin{prop}  \label{prop:bottomdepth2}
Let $\A$ be a depth~$2$, completely decomposable arrangement of $n$ planes, with 
normal form $\A(I_1,\dots ,I_r,J)$.  The bottom characteristic variety $V_{n-2}(\A)$ 
has the following irreducible components:
\begin{alphenum}
\item  Subtori passing through $\mathbf{1}$:  
\begin{align*}
&T({\Check{I}_p}), \ \: \text{for}\ \: p\in [1,r] \quad \text{and}\quad 
T({\Check{J}}), \ \: \text{if}\ \: |J|>1.  
\intertext{
\item
Translated subtori: 
}
&T({\cup_{p\notin P} I_{p}}\cup \{n\}) \cap 
\overline{T}({\cup_{p\in P} \{\max I_p\}}), 
\ \: \text{for}\ \: \emptyset \ne P\subseteq [1,r].
\end{align*}
\end{alphenum}
\end{prop}
 
\begin{exmp}  \label{ex:comps6}
The arrangement $\A=\A(214356)$ from Example~\ref{ex:depth} 
is in normal form, with $I_1=\{2,1\}$, $I_2=\{4,3\}$, and $J=\{5,6\}$.  
The components of $V_4(\A)$ are:
\begin{align*} 
&\{t_6-1=t_4+1=t_2+1=0\},\\
&\{t_6-1=t_4+1=t_2-1=t_1-1=0\},\\
&\{t_6-1=t_4-1=t_3-1=t_2+1=0\},\\
&\{t_6-1=t_5-1=t_4-1=t_3-1=t_2-1=0\},\\
&\{t_6-1=t_5-1=t_4-1=t_2-1=t_1-1=0\},\\
&\{t_6-1=t_4-1=t_3-1=t_2-1=t_1-1=0\}.
\end{align*}
\end{exmp}

\subsection{}
Let $\A$ be a completely decomposable arrangement of 
depth~$2$, with normal form $\A(I_1,\dots ,I_r,J)$.  Recall that 
$|I_1|\le \dots \le |I_r|$.  Define 
\begin{equation*} \label{eq:sofa}
S(\A)=\{|I_1|,\dots ,|I_r|\}.  
\end{equation*}
Since we also have $I_1\le \dots \le I_r<J$, the ordered list $S(\A)$, 
together with the number of planes, $n=|J|+\sum_{k=1}^r |I_k|$, 
determines the normal form.  

Let $\Sigma=\Sigma_{n-2}(\A)$ be the list of codimensions of 
irreducible components of $V_{n-2}(\A)$, and 
$\Sigma_{\mathbf{1}}=\Sigma_{\mathbf{1},n-2}(\A)$ 
be the sublist corresponding to components passing 
through $\mathbf{1}$.   From Proposition~\ref{prop:bottomdepth2}, we see that 
\begin{subequations}  \label{eq:comps}
\begin{alignat}{2}  
\Sigma_{\mathbf{1}}&=
\{n+1-|I_p|\}_{p=1,\dots, r}  && \text{if } |J|=1 \label{eq:c11}\\
\Sigma_{\mathbf{1}}&=
\{n+1-|I_p|\}_{p=1,\dots, r} \cup \{n+1-|J|\} &\qquad& 
\text{if } |J|>1 \label{eq:c12}\\
\Sigma \setminus \Sigma_{\mathbf{1}} &= 
\{r+1+\sum_{p\notin P}(|I_p| -1)\}_{\emptyset\subsetneq P\subseteq [1,r]}  
&\qquad& \label{eq:cminusc1}
\end{alignat}  
\end{subequations}
The lists $\Sigma_{\mathbf{1}}$ and $\Sigma$ have lengths
\begin{equation}  \label{eq:length}
d_1=r+\e_J+1, \qquad d=2^r+r+\e_J,
\end{equation}
where $\e_J=0$ if $|J|>1$ and $\e_J=-1$ if  $|J|=1$. 

\begin{thm}  \label{thm:depth2}
The list $S(\A)$ is a complete homotopy-type invariant 
for depth~$2$, completely decomposable arrangements $\A$ 
of $n$ planes.  
\end{thm}

\begin{proof}  
Let $\A'$ be another completely decomposable arrangement of depth~$2$, 
with normal form  $\A(I'_1, \dots, I'_r, J')$.  Assume $X(\A)\simeq X(\A')$.  
Then, by Theorem~\ref{thm:invariant}, $\Sigma(\A)=\Sigma(\A')$ and 
$\Sigma_{\mathbf{1}}(\A)=\Sigma_{\mathbf{1}}(\A')$.  We want to show that 
$S(\A)=S(\A')$.  There are four cases
to consider, according to the sizes of $J$ and $J'$:

\begin{itemize}
\item $|J|=1, |J'|>1$.  Then, by \eqref{eq:length}, 
the system of equations $d=d'$, $d_1=d_1'$ has no solution.

\item $|J|=|J'|=1$.  Then, by \eqref{eq:c11}, 
$\Sigma_{\mathbf{1}}(\A)=\Sigma_{\mathbf{1}}(\A')$ implies $S(\A)=S(\A')$.  

\item $|J|=|J'|>1$.  Then,  \eqref{eq:c12}, 
$\Sigma_{\mathbf{1}}(\A)=\Sigma_{\mathbf{1}}(\A')$ implies $S(\A)=S(\A')$. 

\item $|J|>|J'|>1$.  Then, by  \eqref{eq:length}, $r=r'$.  If
$r=1$,  equation~\eqref{eq:c12} implies that $\{|I_1|, |J|\}=\{|I'_1|,
|J'|\}$  as {\it unordered} lists.  But condition~\eqref{condition3} from 
Definition~\ref{defn:normalform} and  the fact that $|J|>|J'|$ rule out 
this possibility.  If $r>1$, equation~\eqref{eq:cminusc1} implies 
$|J|=n+r+2-\max( \Sigma \setminus \Sigma_{\mathbf{1}}) - 
\min( \Sigma \setminus \Sigma_{\mathbf{1}})$.  Hence $|J|-|J'|=r-r'=0$, 
which again is impossible.
\end{itemize}

Conversely, assume $S(\A)=S(\A')$.  Then, as noted above, 
the normal forms of $\A$ and $\A'$ coincide.  
By Proposition~\ref{prop:normalform}, $\A$ and $\A'$ are rigidly
isotopic, and thus $X(\A)\simeq X(\A')$.  
\end{proof}

\begin{cor}  \label{cor:depth2}
The number of homotopy classes of $2$-arrangements of $n$~planes 
which are completely decomposable of depth at most ~$2$ equals 
$p(n-1) - \lfloor(n-1)/2\rfloor$, where $p( \cdot)$ is the partition 
function, and $\lfloor\cdot\rfloor$ is the integer part function.  
\end{cor}

\begin{proof}  
Follows from the Theorem by an elementary counting argument.  
\end{proof}


\section{Top Characteristic Varieties} 
\label{sec:topvar}

In this section we study the top characteristic variety $V_1(\A)$ of a 
$2$-arrangement $\A$, and the number $\Tors_{p,1}(\A)$ of its 
$p$-torsion points, for $p$ a prime number.  

\subsection{}  
Let us start with a completely decomposable arrangement, $\A=\A(\t)$.  
Let $\t=\t_0\to\t_1 \to\cdots\to\t_d=(1)$ be the decomposition 
sequence for $\t$.  Recall that each permutation in the sequence
is partitioned into blocks of consecutive integers.  
We call such a block $B$ {\em essential} if either $|B| \ge 2$, or 
$B=\t_{d-1}$ and $|B| > 2$.  

\begin{thm}  \label{thm:topvariety} 
The top characteristic variety $V_1(\A)$ of a completely 
decomposable arrangement of $n$ planes is the union of an arrangement 
$\V(\A)$ of codimension~$1$ subtori of $(\C^*)^n$, all passing 
through $\mathbf{1}$.  
\end{thm}
\begin{proof} 
Choose $\t\in S_n$ so that $\A=\A(\t)$ and $\depth(\A)=d(\t)$.  
Let $m$ be the number of essential blocks of $\t$.   
Recall from \S\ref{subsec:shift} that $L=L(\A)$ is an iterated torus link, 
obtained by $(1,\pm 1)$-cablings on the unknot.  Thus, it is 
a spliced link in the sense of Eisenbud and Neumann~\cite{EN}.  
The decomposition sequence of $\t$ corresponds to a minimal 
splice diagram of $L$:  The (signed) essential blocks $B_{1},\dots , B_{m}$ 
of the permutations in the sequence correspond to the (signed) nodes 
$v_{n+1},\dots , v_{n+m}$ of the diagram, and the integers $1, \dots, n$ 
to the arrowheads $v_{1},\dots , v_{n}$.   Then, according to 
\cite{EN}, Theorem~12.1, the Alexander polynomial of $L$ is given by:
\[
\D_{L}(t_1, \dots, t_n)=\prod_{j=n+1}^{n+m}(t_1^{l_{1,j}}t_2^{l_{2,j}} 
\cdots t_n^{l_{n,j}} -1)^{\delta_j -2},
\]
where $l_{i,j}=\pm 1$ is the linking number of $L_i$ with the 
``virtual component" corresponding to $v_j$, and $\delta_j$ 
is the valency of $v_j$.  Thus, each irreducible component of 
$V_1(\A)=\{\D_{L}=0\}$ is a codimension~$1$ subtorus.  
(It actually can be shown that $V_1(\A)$ has precisely $m$ components.)
\end{proof} 

To compute the number of torsion points on $V_1(\A)$, we may 
now use a result of Bj\"{o}rner and Ekedahl~\cite{BE}.  
Indeed, an arrangement $\V$ of codimension~$1$ subtori in 
$(\C^*)^n$ defines an arrangement $\V_{p}$ of hyperplanes 
in $(\Z_p)^n$:  
To a subtorus $t_1^{a_1}\cdots t_n^{a_n}-1=0$ corresponds 
the hyperplane $a_1 x_1+\cdots +a_nx_n=0 \mod p$.   
Proposition~3.2 of \cite{BE} then implies the following.  

\begin{prop}	  \label{prop:p-torsion}  
The number of $p$-torsion points on the union 
$U=U(\V)$ of an arrangement of subtori in $(\C^*)^n$ is given by:
\[
 \Tors_p(U)=-\sum_{x \in L \setminus \{\hat{0}\}}~\mu(\hat{0},x)p^{\dim(x)}, 
\]
where $L$ is the intersection lattice of the arrangement 
$\V_{p}$, with minimal element $\hat{0}=(\Z_p)^n$, and 
M\"{o}bius function $\mu$. 
\end{prop}

\begin{exmp} \label{ex:depth3}
The arrangement $\A=\A(312546)$ is completely decomposable, of depth $3$.
Its decomposition sequence is $(312546)\to (2134)\to (12)\to (1)$.  
Proposition~\ref{prop:combed} gives $\xi=A_{1,3}A_{2,3}A_{4,5}$. 
The Artin representation of $\xi$ in the basis 
$y_1=x_1$, $y_2=x_1x_2$, $y_3=x_1x_2x_3$, $y_4=x_4$, $y_5=x_4x_5$ is 
given by $\xi(y_1)=y_3y_2^{-1}y_1y_2y_3^{-1}$, $\xi(y_2)=y_3y_2y_3^{-1}$,
$\xi(y_3)=y_3$, $\xi(y_4)=y_5y_4y_5^{-1}$, and $\xi(y_5)=y_5$.
The Alexander polynomial is $\D_{\A}(t_1,\dots,t_6)=
(t_6-1)(t_6-t_5^{2})(t_6-t_3^{2})(t_6-t_3^{2}t_2^{-2})$.  
Proposition~\ref{prop:p-torsion} yields $\Tors_{2,1}(\A)=32$ 
and $\Tors_{3,1}(\A)=585$.
\end{exmp}

For an arrangement of depth~$2$, we can give a more precise description 
of the top characteristic variety, based on formula \eqref{eq:alexdepth2} 
for the Alexander matrix.
\begin{prop}  \label{prop:torsdecomp}
Let $\A$ be a depth~$2$ arrangement of $n$ planes, 
with normal form $\A(I_1, \dots, I_r, J)$, and let 
$k_q=\max I_q$, for $1 \le q \le r$.  Then: 
\vskip 5pt
\begin{alphenum} 
\item 
$\Delta_{\A}(t_1,\dots,t_n)=
(t_n-1)^{|J|+r-2}\prod_{q=1}^{r}(t_n-t_{k_q}^2)^{|I_q|-1};$
\vskip 4pt
\item
$V_{1}(\A)=\{\mathbf{t}\in(\C^*)^n \mid 
(t_n-1)\prod_{q=1}^{r}(t_n-t_{k_q}^2) = 0\};$
\vskip 4pt
\item
$\Tors_{2,1}(\A)=2^{n-1}$ and 
$\Tors_{p,1}(\A)=p^{n-r-1}(p^{r+1}-(p-1)^{r+1})$, for $p\ge3$.
\end{alphenum}
\end{prop}

\subsection{}
Let $L$ be a link in $\SP^3$.  The Alexander polynomial of a sublink 
of $L$, and that of an $(a,b)$-cable about $L$, can be computed from 
the Alexander polynomial of $L$, via the following well-known formulae of 
Torres and Sumners--Woods, see~\cite{EN, Hi, T1}. 

\begin{thm} \label{thm:tsw}
Let $L=L_1 \cup \dots \cup L_n$ be a link in $\SP^3$.  
Set $T=t_1^{l_1} \cdots t_{n-1}^{l_{n-1}}$, where $l_i=\lk(L_i,L_n)$.  
Then:
\begin{subequations} \label{eq:iterated}
\begin{align}  
\D_L(t_{1}, \dots , t_{n-1}, 1)&=(T-1)\D_{L \setminus L_n}
(t_{1}, \dots , t_{n-1}), 
\label{eq:torres}\\
\intertext{
Moreover, if $L'=L\{a,b\}$, with $\gcd(a,b)=1$, then:}
\D_{L'}(t_{1}, \dots , t_{n}, t_{n+1})&=
(T^{a}t_{n}^{b}t_{n+1}^{b}-1)\D_{L}(t_{1}, \dots , t_{n-1}, t_n^{a}t_{n+1}).
 \label{eq:sw} 
\end{align}
\end{subequations}
\end{thm}

\begin{cor}  \label{cor:vcable}
Let $\A$ be an arrangement of $n$ planes, and let $\A^k\{r\}$ 
be an $r$-cable about it.  Then:
\vskip 4pt
\begin{alphenum} 
\item  \label{it1}
$V_{1}(\A\{r\})=\{\bt\in 
(\C^*)^{n+r}\mid t_1\cdots t_n -1=0 \,\text{ or }\, 
(t_1,\dots ,t_n)\in V_{1}(\A)\};$
\vskip 4pt
\item  \label{it2}
$\Tors_{p,1}(\A\{r\})=p^{r-1}\Tors_{p,1}(\A\{1\}).$
\end{alphenum}
\end{cor}
\begin{proof}  
Let $L=L(\A)$.  
Recall from \S\ref{cables} that $L(\A\{1\})=L\{1,1\}$.  
By the Sumners-Woods formula \eqref{eq:sw}, we have 
\[
\D_{L(\A\{1\})}(t_1, \dots , t_n, t_{n+1})=
(t_1\cdots t_{n+1}-1)\D_{L}(t_1, \dots , t_{n-1}, t_nt_{n+1}). 
\]
After a monomial change of basis, 
this implies part \eqref{it1} for $r=1$.  
The general case follows from the same formula, by induction on $r$.  
Part \eqref{it2} follows immediately from part \eqref{it1}.
\end{proof}

\subsection{}
We conclude this section with a recursion formula for $\Tors_{2,1}(\A)$.  
Let $\Delta_{L(\A)}$ be the Alexander polynomial of the link 
$L(\A)$.  Define the {\em single-variable Alexander polynomial} of $\A$ to be 
$\Delta_{\A}(t):=(t-1)\Delta_{L(\A)}(t,\dots ,t)$.  Furthermore, set 
\[
\delta(\A)=\begin{cases}
1&\quad\text{if } \D_{\A}(-1)=0, \\
0&\quad\text{otherwise}. 
\end{cases}
\]
\begin{exmp}  \label{ex:deltacx}
The single-variable Alexander polynomial of the complex arrangement $\A_n$ is 
$\D_{\A_n}(t)=(t-1)(t^n-1)^{n-2}$.  Thus, $\delta(\A_n)=\frac{1+(-1)^n}{2}$.  
\end{exmp}

\begin{thm} \label{thm:2tors}
Let $\A=\{H_1,\dots, H_n\}$ be a $2$-arrangement in $\R^4$.   
The number of $2$-torsion points of the top characteristic 
variety $V_1(\A)$ is given by the following formula:
\begin{equation} \label{eq:deltasum}
\Tors_{2,1}(\A)=2^{n-1}-\frac{1+(-1)^n}{2}+\delta(\A)+
\sum_{\B\in \Gamma(\A)}\delta(\B),
\end{equation}
where $\Gamma(\A)$ is the set of all indecomposable, proper 
sub-arrangements of $\A$ with an odd number of planes.
\end{thm}
\begin{proof}  
 By definition, 
\begin{equation} \label{eq:torsum}
\Tors_{2,1}(\A) = \sum_{\omega\in\Omega_2^n} c_{\omega},
\end{equation}
where $\Omega_2^n = \{(\omega_1,\dots ,\omega_n)\in(\C^*)^n \mid \omega_i=\pm 1\}$,  
and $c_{\omega}=1$ if $\Delta_{L(\A)}(\omega)=0$,   
and $c_{\omega}=0$ otherwise.  For $\omega\in\Omega_2^n$, let 
$\Az=\{ H_i\in\A \mid \omega_i = -1\}$.  
There are several cases to consider:

\begin{itemize}
\item If $\omega=(-1,\dots ,-1)$, then $\Az=\A$, and so 
$c_{\omega}=\delta(\A)$.

\item If $\omega\ne (-1,\dots ,-1)$, then $\Az$ is a proper 
sub-arrangement of $\A$, and so, by repeated application of 
Torres's formula~\eqref{eq:torres}, we have
\begin{equation} \label{eq:delomega}
\Delta_{L(\A)}(\omega) = ((-1)^{|\Az|} - 1)\Delta_{L(\Az)}(-1,\dots ,-1).
\end{equation}
\begin{itemize}
\item
If $|\Az|$ is even, this formula says that 
$\Delta_{L(\A)}(\omega)=0$, and so $c_{\omega}=1$.  
There are $2^{n-1}-\frac{1+(-1)^n}{2}$ such contributions to 
the sum \eqref{eq:torsum}.    

\item
If $|\Az|$ is odd, and $\Az$ is decomposable, we may write  
$\Az=\A_{\omega}'\{\pm 1\}$, with the cabling done about the last component 
of $\A_{\omega}'$.  Let $\A_{\omega}''$ be the sub-arrangement obtained by deleting 
the last component of $\A_{\omega}'$.  Clearly, $|\A_{\omega}''|=|\Az|-2$.   
Formulas \eqref{eq:sw}, and \eqref{eq:torres} give
\begin{equation}  \label{eq:recursion}
\Delta_{L(\Az)}(-1,\dots ,-1)=-2\Delta_{L(\A_{\omega}')}(-1,\dots ,-1,1) 
=4\Delta_{L(\A_{\omega}'')}(-1,\dots ,-1).
\end{equation}
Hence, $\delta(\Az)=\delta(\A_{\omega}'')$.   Iterating this decabling-deletion 
procedure, we eventually reach an arrangement $\B$ for which the procedure 
must stop.  There are two possibilities:  
\begin{itemize}
\item 
One is $\B=\A_3$, in which case $c_{\omega}=\delta(\A_3)=0$.   

\item  
The other is $\B\in \Gamma(\A)$, in which case $c_{\omega}=\delta(\B)$.  
Clearly, any element of $\Gamma(\A)$ can be reached by the above 
procedure; thus there are $|\Gamma(\A)|$ such contributions to the sum 
\eqref{eq:torsum}. 
\end{itemize} 
\end{itemize}
\end{itemize}
This completes the proof.  
\end{proof}

\begin{rem}  \label{rem:lowboundtors}
Note that $ 2^{n-1} -1\le \Tors_{2,1}(\A)\le 2^{n}$.  If $\Tors_{2,1}(\A) = 2^{n-1}-1$, 
and $n\ge 3$, then the top characteristic variety $V_1(\A)$ is {\em not} 
the union of translated subtori of $(\C^*)^n$.   For, otherwise, 
at least one of the subtori must be of the form $T=\{t_1^{a_1}\cdots t_n^{a_n}-1=0\}$, 
since $\D_{L(\A)}(1,\dots,1)=0$.  But the torus $T$ has $2^{n-1}$ torsion 
points of order $2$.  
\end{rem}

\begin{cor} \label{cor:almostcompldec}
If all the proper subarrangements of $\A$ are completely decomposable, 
then $\Tors_{2,1}(\A)=2^{n-1}-\frac{1+(-1)^n}{2}+\delta(\A)$.  
\end{cor}

\begin{cor} \label{cor:compldec}
If $\A$ is completely decomposable, then $\Tors_{2,1}(\A)=2^{n-1}$. 
\end{cor}
\begin{proof}
The recursion formula~\eqref{eq:recursion}, together with 
Example~\ref{ex:deltacx} imply that $\delta(\A)=\frac{1+(-1)^n}{2}$, 
and the conclusion follows from the previous corollary. 
\end{proof}

\begin{exmp} \label{ex:h1}
The arrangement $\A=\A(31425)$ is horizontal, indecomposable, and all 
its subarrangements are completely decomposable.  The (single variable) 
Alexander polynomial is $\D_{\A}(t)=(t-1)^4(4t^2-t+4)$, and so $\delta(\A)=0$.  
From Corollary~\ref{cor:almostcompldec}, we get $\Tors_{2,1}(\A)=16$.  
\end{exmp}

\begin{exmp} \label{ex:h2}
The arrangement $\A=\A(314256)$ is decomposable, but not completely 
decomposable, since it has $\A(31425)$ as a subarrangement. We have 
$\D_{\A}(t)=(t^6-1)(t-1)^3(t+1)(3t^2-2t+3)$, and so $\delta(\A)=1$.  
From Theorem~\ref{thm:2tors}, we get $\Tors_{2,1}(\A)=32$.  
\end{exmp}

\begin{exmp} \label{ex:h3}
The arrangement $\A=\A(241536)$ is horizontal, indecomposable, and all 
its proper subarrangements are completely decomposable.  We have 
$\D_{\A}(t)=(t-1)^5(5t^4+6t^2+5)$, and so $\delta(\A)=0$.  
From Corollary~\ref{cor:almostcompldec}, we get $\Tors_{2,1}(\A)=31$. 
Hence, $V_1(\A)$ is not a union of translated subtori of $(\C^*)^6$.  
\end{exmp}


\section{Mazurovski\u{\i}'s arrangements} 
\label{sec:maz}

In this section, we study the $2$-arrangements associated to 
Mazurovski\u{\i}'s configurations.  Using their associated cablings,   
we find infinitely many pairs of arrangements whose complements are 
cohomologically equivalent, but not homotopy equivalent.

\subsection{}  
In \cite{M1}, Mazurovski\u{\i} introduced a remarkable pair of 
configurations of skew lines, $K$ and $L$, which have the same 
linking numbers, but are not rigidly isotopic.  Let $\KK=\A(K)$ and 
$\LL=\A(L)$ be the corresponding arrangements of planes.   
The arrangement $\KK$ is horizontal, with associated permutation 
$\tau=(341256)$.  Moreover, $\KK$ is completely decomposable, of depth $3$; 
a minimal decomposition sequence is $(341256)\to (213)\to (12)\to (1)$.
The arrangement $\LL$ is neither horizontal, nor decomposable.  
Defining polynomials for $\KK$ and $\LL$ are given by 
\begin{align*}  \label{polyKL}
f_{\KK}(z,w)&=f(z,w)\cdot (z-7w),\\
f_{\LL}(z,w)&=f(z,w)\cdot (z-\frac{6-7\ii}{2} w-
\frac{3+14\ii}{2}\bar{w}),
\end{align*}
where $f$ is the following defining polynomial for $\A(34125)$:   
\begin{equation*}  \label{polyf}
\begin{split}
f(z,w)&=(z-\frac{5-5\ii}{2}w+\frac{3-5\ii}{2}\bar{w})
(z-\frac{7-10\ii}{2}w+\frac{3-10\ii}{2}\bar{w})\\
&\quad\times
(z-\frac{5-14\ii}{2}w-\frac{3+14\ii}{2}\bar{w})
(z-\frac{7-9\ii}{2}w-\frac{3+9\ii}{2}\bar{w})
(z-6w).\\
\end{split}
\end{equation*}

\begin{figure}
\[
\beginsmallgraph
\Label {$1$} at (1, 8)
\Label {$2$} at (1, 9)
\Label {$3$} at (1, 10)
\Label {$4$} at (1, 11)
\Label {$5$} at (1, 12)
\Label {$6$} at (1, 13)
\edge from (2, 8) to (3, 8)
\edge from (2, 9) to (3, 9)
\edge from (2, 10) to (3, 10)
\edge from (2, 11) to (3, 11)
\edge from (2, 12) to (3, 12)
\edge from (2, 13) to (3, 13)
\edge from (3, 8) to (4, 8)
\edge from (3, 9) to (4, 9)
\edge from (3, 10) to (4, 10)
\edge from (3, 11) to (4, 11)
\edge from (3, 12) to (3.3, 12.3)
\edge from (3.7, 12.7) to (4, 13)
\edge from (3, 13) to (4, 12)
\edge from (4, 8) to (5, 8)
\edge from (4, 9) to (5, 9)
\edge from (4, 10) to (5, 10)
\edge from (4, 11) to (4.3, 11.3)
\edge from (4.7, 11.7) to (5, 12)
\edge from (4, 12) to (5, 11)
\edge from (4, 13) to (5, 13)
\edge from (5, 8) to (6, 8)
\edge from (5, 9) to (6, 9)
\edge from (5, 10) to (5.3, 10.3)
\edge from (5.7, 10.7) to (6, 11)
\edge from (5, 11) to (6, 10)
\edge from (5, 12) to (6, 12)
\edge from (5, 13) to (6, 13)
\edge from (6, 8) to (7, 8)
\edge from (6, 9) to (6.3, 9.3)
\edge from (6.7, 9.7) to (7, 10)
\edge from (6, 10) to (7, 9)
\edge from (6, 11) to (7, 11)
\edge from (6, 12) to (7, 12)
\edge from (6, 13) to (7, 13)
\edge from (7, 8) to (7.3, 8.3)
\edge from (7.7, 8.7) to (8, 9)
\edge from (7, 9) to (8, 8)
\edge from (7, 10) to (8, 10)
\edge from (7, 11) to (8, 11)
\edge from (7, 12) to (8, 12)
\edge from (7, 13) to (8, 13)
\edge from (8, 8) to (9, 8)
\edge from (8, 9) to (9, 9)
\edge from (8, 10) to (9, 10)
\edge from (8, 11) to (9, 11)
\edge from (8, 12) to (8.3, 12.3)
\edge from (8.7, 12.7) to (9, 13)
\edge from (8, 13) to (9, 12)
\edge from (9, 8) to (10, 8)
\edge from (9, 9) to (10, 9)
\edge from (9, 10) to (10, 10)
\edge from (9, 11) to (9.3, 11.3)
\edge from (9.7, 11.7) to (10, 12)
\edge from (9, 12) to (10, 11)
\edge from (9, 13) to (10, 13)
\edge from (10, 8) to (11, 8)
\edge from (10, 9) to (11, 9)
\edge from (10, 10) to (10.3, 10.3)
\edge from (10.7, 10.7) to (11, 11)
\edge from (10, 11) to (11, 10)
\edge from (10, 12) to (11, 12)
\edge from (10, 13) to (11, 13)
\edge from (11, 8) to (12, 8)
\edge from (11, 9) to (11.3, 9.3)
\edge from (11.7, 9.7) to (12, 10)
\edge from (11, 10) to (12, 9)
\edge from (11, 11) to (12, 11)
\edge from (11, 12) to (12, 12)
\edge from (11, 13) to (12, 13)
\edge from (12, 8) to (13, 8)
\edge from (12, 9) to (13, 9)
\edge from (12, 10) to (13, 10)
\edge from (12, 11) to (13, 11)
\edge from (12, 12) to (12.3, 12.3)
\edge from (12.7, 12.7) to (13, 13)
\edge from (12, 13) to (13, 12)
\edge from (13, 8) to (14, 8)
\edge from (13, 9) to (14, 9)
\edge from (13, 10) to (14, 10)
\edge from (13, 11) to (14, 12)
\edge from (13, 12) to (13.3, 11.7)
\edge from (13.7, 11.3) to (14, 11)
\edge from (13, 13) to (14, 13)
\edge from (14, 8) to (15, 8)
\edge from (14, 9) to (15, 9)
\edge from (14, 10) to (15, 11)
\edge from (14, 11) to (14.3, 10.7)
\edge from (14.7, 10.3) to (15, 10)
\edge from (14, 12) to (15, 12)
\edge from (14, 13) to (15, 13)
\edge from (15, 8) to (16, 8)
\edge from (15, 9) to (16, 9)
\edge from (15, 10) to (16, 10)
\edge from (15, 11) to (16, 11)
\edge from (15, 12) to (16, 13)
\edge from (15, 13) to (15.3, 12.7)
\edge from (15.7, 12.3) to (16, 12)
\edge from (16, 8) to (17, 8)
\edge from (16, 9) to (17, 9)
\edge from (16, 10) to (17, 10)
\edge from (16, 11) to (17, 12)
\edge from (16, 12) to (16.3, 11.7)
\edge from (16.7, 11.3) to (17, 11)
\edge from (16, 13) to (17, 13)
\edge from (17, 8) to (18, 8)
\edge from (17, 9) to (18, 9)
\edge from (17, 10) to (18, 10)
\edge from (17, 11) to (18, 11)
\edge from (17, 12) to (17.3, 12.3)
\edge from (17.7, 12.7) to (18, 13)
\edge from (17, 13) to (18, 12)
\edge from (18, 8) to (19, 8)
\edge from (18, 9) to (19, 9)
\edge from (18, 10) to (19, 10)
\edge from (18, 11) to (19, 11)
\edge from (18, 12) to (19, 12)
\edge from (18, 13) to (19, 13)
%
%
\Label {$1$} at (1, 1)
\Label {$2$} at (1, 2)
\Label {$3$} at (1, 3)
\Label {$4$} at (1, 4)
\Label {$5$} at (1, 5)
\Label {$6$} at (1, 6)
\edge from (2, 1) to (3, 1)
\edge from (2, 2) to (3, 2)
\edge from (2, 3) to (3, 3)
\edge from (2, 4) to (3, 4)
\edge from (2, 5) to (3, 5)
\edge from (2, 6) to (3, 6)
\edge from (3, 1) to (4, 1)
\edge from (3, 2) to (4, 2)
\edge from (3, 3) to (4, 3)
\edge from (3, 4) to (4, 4)
\edge from (3, 5) to (3.3, 5.3)
\edge from (3.7, 5.7) to (4, 6)
\edge from (3, 6) to (4, 5)
\edge from (4, 1) to (5, 1)
\edge from (4, 2) to (5, 2)
\edge from (4, 3) to (5, 3)
\edge from (4, 4) to (4.3, 4.3)
\edge from (4.7, 4.7) to (5, 5)
\edge from (4, 5) to (5, 4)
\edge from (4, 6) to (5, 6)
\edge from (5, 1) to (6, 1)
\edge from (5, 2) to (6, 2)
\edge from (5, 3) to (5.3, 3.3)
\edge from (5.7, 3.7) to (6, 4)
\edge from (5, 4) to (6, 3)
\edge from (5, 5) to (6, 5)
\edge from (5, 6) to (6, 6)
\edge from (6, 1) to (7, 1)
\edge from (6, 2) to (6.3, 2.3)
\edge from (6.7, 2.7) to (7, 3)
\edge from (6, 3) to (7, 2)
\edge from (6, 4) to (7, 4)
\edge from (6, 5) to (7, 5)
\edge from (6, 6) to (7, 6)
\edge from (7, 1) to (7.3, 1.3)
\edge from (7.7, 1.7) to (8, 2)
\edge from (7, 2) to (8, 1)
\edge from (7, 3) to (8, 3)
\edge from (7, 4) to (8, 4)
\edge from (7, 5) to (8, 5)
\edge from (7, 6) to (8, 6)
\edge from (8, 1) to (9, 1)
\edge from (8, 2) to (9, 2)
\edge from (8, 3) to (9, 3)
\edge from (8, 4) to (9, 4)
\edge from (8, 5) to (8.3, 5.3)
\edge from (8.7, 5.7) to (9, 6)
\edge from (8, 6) to (9, 5)
\edge from (9, 1) to (10, 1)
\edge from (9, 2) to (10, 2)
\edge from (9, 3) to (9.3, 3.3)
\edge from (9.7, 3.7) to (10, 4)
\edge from (9, 4) to (10, 3)
\edge from (9, 5) to (10, 5)
\edge from (9, 6) to (10, 6)
\edge from (10, 1) to (11, 1)
\edge from (10, 2) to (11, 2)
\edge from (10, 3) to (11, 3)
\edge from (10, 4) to (11, 5)
\edge from (10, 5) to (10.3, 4.7)
\edge from (10.7, 4.3) to (11, 4)
\edge from (10, 6) to (11, 6)
\edge from (11, 1) to (12, 1)
\edge from (11, 2) to (12, 2)
\edge from (11, 3) to (11.3, 3.3)
\edge from (11.7, 3.7) to (12, 4)
\edge from (11, 4) to (12, 3)
\edge from (11, 5) to (12, 5)
\edge from (11, 6) to (12, 6)
\edge from (12, 1) to (13, 1)
\edge from (12, 2) to (13, 3)
\edge from (12, 3) to (12.3, 2.7)
\edge from (12.7, 2.3) to (13, 2)
\edge from (12, 4) to (13, 4)
\edge from (12, 5) to (13, 5)
\edge from (12, 6) to (13, 6)
\edge from (13, 1) to (14, 1)
\edge from (13, 2) to (14, 2)
\edge from (13, 3) to (13.3, 3.3)
\edge from (13.7, 3.7) to (14, 4)
\edge from (13, 4) to (14, 3)
\edge from (13, 5) to (14, 5)
\edge from (13, 6) to (14, 6)
\edge from (14, 1) to (15, 1)
\edge from (14, 2) to (15, 2)
\edge from (14, 3) to (15, 3)
\edge from (14, 4) to (15, 4)
\edge from (14, 5) to (15, 6)
\edge from (14, 6) to (14.3, 5.7)
\edge from (14.7, 5.3) to (15, 5)
\edge from (15, 1) to (16, 1)
\edge from (15, 2) to (16, 2)
\edge from (15, 3) to (16, 3)
\edge from (15, 4) to (16, 5)
\edge from (15, 5) to (15.3, 4.7)
\edge from (15.7, 4.3) to (16, 4)
\edge from (15, 6) to (16, 6)
\edge from (16, 1) to (17, 1)
\edge from (16, 2) to (17, 2)
\edge from (16, 3) to (16.3, 3.3)
\edge from (16.7, 3.7) to (17, 4)
\edge from (16, 4) to (17, 3)
\edge from (16, 5) to (17, 5)
\edge from (16, 6) to (17, 6)
\edge from (17, 1) to (18, 1)
\edge from (17, 2) to (18, 2)
\edge from (17, 3) to (18, 3)
\edge from (17, 4) to (18, 4)
\edge from (17, 5) to (17.3, 5.3)
\edge from (17.7, 5.7) to (18, 6)
\edge from (17, 6) to (18, 5)
\edge from (18, 1) to (19, 1)
\edge from (18, 2) to (19, 2)
\edge from (18, 3) to (19, 3)
\edge from (18, 4) to (19, 4)
\edge from (18, 5) to (19, 5)
\edge from (18, 6) to (19, 6)
\endgraph
\]
\caption{Mazurovski\u\i's pair: The half-braids $\a_{\KK}$ (top) and 
$\a_{\LL}$ (bottom).} 
\label{fig:KandL}
\end{figure}

The half-braids associated to $\KK$ and $\LL$ are pictured in 
Figure~\ref{fig:KandL}.  We see that the linking numbers of $L(\KK)$ are
$l_{1,4}=l_{2,4}=l_{1,3}=l_{2,3}=-1$, and all other $l_{i,j}=1$, 
whereas the linking numbers of $L(\LL)$ are 
$l_{1,5}=l_{2,5}=l_{1,4}=l_{2,4}=-1$, 
and all other $l_{i,j}=1$.  The reordering of the components of $L(\LL)$ 
that fixes $1$, $2$, $6$ and permutes $3$, $4$, $5$ to $4$, $5$, $3$ 
identifies the linking numbers of $L(\KK)$ and $L(\LL)$.   Thus, 
$H^*(X(\KK);\Z)\cong H^*(X(\LL);\Z)$.

\subsection{} \label{subsection:kandl}
In order to distinguish between the cohomologically equivalent 
arrangements $\KK$ and $\LL$, we turn to their characteristic varieties.   
From Figure~\ref{fig:KandL}, we see that the reduced half-braids 
of $\KK$ and $\LL$ are:
\begin{equation*}  \label{halfKL}
\check{\a}_{\KK}=\s_4\s_3\s_2\s_1\s_4\s_3^{-1}\s_2^{-1}\s_4^{-1}\s_3^{-1}\s_4, \quad
\check{\a}_{\LL}=\s_4\s_2\s_3^{-1}\s_2\s_1^{-1}\s_2\s_4^{-1}\s_3^{-1}\s_2\s_4.
\end{equation*}
The braids $\xi=\Delta_5 \check{\a}^{-1}\in P_5$ are expressed in terms 
of the pure braid generators, as follows.  For $\KK$, which is horizontal, 
Proposition~\ref{prop:combed} yields $\xi_{\KK}=A_{1,3}A_{2,3}A_{1,4}A_{2,4}$. 
For $\LL$, it is more convenient to work with the conjugate 
$\xi'_{\LL}=\delta^{-1} \xi_{\LL}\delta $, where $\delta=\s_1\s_3\s_4^{-1}$.  
Routine combing of the braid yields 
$\xi'_{\LL}=A_{1,3}A_{2,3}A_{4,5}A_{1,4}A_{4,5}^{-1}A_{2,4}$.

The Artin representation of $\xi=\xi_{\KK}$ in the basis 
$y_1=x_1$, $y_2=x_3$, $y_3=x_1x_2$, $y_4=x_1x_2x_3x_4$, $y_5=x_5$ is 
given by $\xi(y_1)=y_4y_3^{-1}y_1y_3y_4^{-1}$, $\xi(y_2)=y_3y_2y_3^{-1}$,
$\xi(y_3)=y_4y_3y_4^{-1}$, $\xi(y_4)=y_4$, and $\xi(y_5)=y_5$.
The Alexander matrix of $\KK$ is:
\begin{equation*}  \label{eq:alexK}
\begin{pmatrix}
t_6-t_4^{2}t_3^{-2}  & 0 &  t_4^{2}t_3^{-2}(t_3+1)(1-t_1) 	& (t_4+1)(t_1-1) & 0 & t_1-1 \\ 
0             							& t_6-t_3^{2} & (t_4+t_3)(t_2-1)						& (1-t_3)(t_2-1)	& 0 & t_2-1 \\
0             							& 0 & t_6-t_4^{2}  																			& (t_4+1)(t_3-1) & 0 & t_3-1 \\ 
0             							& 0 & 0                 														&					t_6-1						& 0 & t_4-1 \\
0             							& 0 & 0                 														& 0										& t_6-1 & t_5-1 \\
\end{pmatrix}
\end{equation*}

An elementary computation shows that the bottom variety $V_4(\KK)$ 
has $6$ irreducible components---$3$ codimension~$4$ translated 
subtori of $(\C^{*})^6$, and $3$ codimension~$5$ subtori passing 
through $\mathbf{1}$---given by the following equations:
\begin{align*} 
&\{t_6-1=t_4+1=t_3+1=t_2-1=0\},\\
&\{t_6-1=t_4+1=t_3-1=t_1-1=0\},\\
&\{t_6-1=t_5-1=t_4-1=t_3+1=0\},\\
&\{t_6-1=t_5-1=t_4-1=t_3-1=t_1-1=0\},\\
&\{t_6-1=t_5-1=t_4-1=t_3-1=t_2-1=0\},\\
&\{t_6-1=t_4-1=t_3-1=t_2-1=t_1-1=0\}.
\end{align*} 

The primary decomposition of the ideal $E_4(\LL)$ is much harder to find.  
The implementation in {\it Macaulay~2}~\cite{GS} of the Eisenbud, Huneke, 
and Vasconcelos algorithm yields such a decomposition, and the result is 
that $V_4(\LL)=V_4(\KK)$.  Thus, the bottom varieties fail to distinguish 
between the $\KK$ and $\LL$ arrangements.  

Let us then consider the top varieties.  It is readily seen that 
the Alexander polynomial of $\KK$ is $\Delta_{\KK}(t_1,\dots,t_6)=
(t_6-1)(t_6-t_3^{2})(t_6-t_4^{2})(t_6-t_4^{2}t_3^{-2})$, and so 
$V_1(\KK)$ is the union of $4$ codimension~$1$ subtori 
of $(\C^{*})^6$.  Since $\KK$ is completely decomposable, 
Corollary~\ref{cor:compldec} implies that $\Tors_{2,1}(\KK)=32$.  
 
The Alexander polynomial of $\LL$ may be computed using 
{\it Mathematica}~\cite{W}.  The result is too long to be displayed here, 
but suffices to say that it is an irreducible polynomial over $\Z$, 
consisting of $667$ monomials.   Direct computation shows that the 
single variable Alexander polynomial is 
$\Delta_{\LL}(t)=3(t-1)^5(3t^2 - 2t + 3)^2$.  Hence $\delta({\LL})=0$. 
Since, as is readily checked, all proper subarrangements of $\LL$ are
completely decomposable, Corollary~\ref{cor:almostcompldec}
implies that $\Tors_{2,1}(\LL)=31$. 

Thus, $\Tors_{2,1}(\KK)\ne\Tors_{2,1}(\LL)$.  (As noted in 
Remark~\ref{rem:lowboundtors}, this arithmetic difference 
translates into a geometric  difference: $V_1(\KK)$ is a 
union of subtori, whereas $V_1(\LL)$ is not  even the union 
of translated subtori.)   Appealing now to Theorem~\ref{thm:invariant}, 
we conclude that the complements of $\KK$ and $\LL$ are not homotopy 
equivalent, although, as mentioned previously, they are cohomologically 
isomorphic.  This answers Ziegler's question from \cite{Z}.

\subsection{}
We now use cablings of $\KK$ and $\LL$ to show that the above phenomenon 
happens for arrangements of $n$ planes, for any $n\ge 6$.

\begin{thm}  \label{thm:mazcables}
Let $\KK$ and $\LL$ be Mazurovski\u{\i}'s arrangements of $6$ 
transverse planes in $\R^4$. Let $\KK\{r\}$ and $\LL\{r\}$ be their $r$-cables. 
Then, for each $r\ge 0$, 
\begin{alphenum}
\item 
$H^*(X(\KK\{r\});\Z)\cong H^*(X(\LL\{r\});\Z)$;
\vskip 4pt
\item 
$X(\KK\{r\})\not\simeq X(\LL\{r\})$.
\end{alphenum}
\end{thm}

\begin{proof}  As noted above, the links of $\KK$ and $\LL$ 
have the same linking numbers.  Hence, the links of $\KK\{r\}$ and 
$\LL\{r\}$ have the same linking numbers.  This implies that the complements 
of $\KK\{r\}$ and $\LL\{r\}$ are cohomologically equivalent.

Although $\KK$ and $\LL$ are distinguished by their $2$-torsion points, 
their cables are not.  Indeed, for $r \ge 1$, 
$\Tors_{2,1}(\KK\{r\})=\Tors_{2,1}(\LL\{r\})=2^{r+5}$, 
as can be deduced from Theorem~\ref{thm:2tors} for $r=1$, 
and from Corollary~\ref{cor:vcable} for $r>1$.   
Hence, we turn to $3$-torsion points.  A {\it Mathematica} computation 
shows that $\Tors_{3,1}(\mathcal K\{1\}) = 3^5\cdot 7$ and 
$\Tors_{3,1}(\mathcal L\{1\}) = 3^3\cdot 61$.  
From Corollary~\ref{cor:vcable}, we get 
\[
\Tors_{3,1}(\mathcal K\{r\}) = 3^{r+4}\cdot 7,  
\quad\text{and}\quad 
\Tors_{3,1}(\mathcal L\{r\}) =  3^{r+2}\cdot 61, 
\]
showing that the respective complements are indeed 
not homotopy equivalent. 
\end{proof}

\subsection{}
Mazurovski\u{\i} introduced in \cite{M1} another interesting configuration 
of $6$ lines, which he called $M$.  Like the $L$ configuration, 
the $M$ configuration is non-horizontal and indecomposable (they  
are the only two such configurations of $6$ lines, up to rigid isotopy 
and mirror images).  But, unlike $L$, the $M$ configuration does 
not have the linking numbers of any horizontal configuration.  
Let $\MM=\A(M)$ be the corresponding arrangement.  
A defining polynomial for it is:
\begin{equation*}  \label{Mpoly}
\begin{split}
f_{\MM}(z,w)&=
(z-(10-\ii)w+
(9-4\ii)\bar{w})
(z-(3-4\ii)w-
(1+4\ii)\bar{w})
\\
&\quad\times
(z-\frac{5-10\ii}{2}w+
\frac{1-10\ii}{2}\bar{w})
(z-(6-5\ii)w+
\frac{1-10\ii}{2}\bar{w})
\\
&\quad\times
(z-\frac{21-29\ii}{4}w+
\frac{1-9\ii}{4}\bar{w})
(z-6w).
\end{split}
\end{equation*}
\noindent
The reduced braids associated to $\MM$ are:
\[
\check{\a}_{\MM}=\s_2\s_3^{-1}\s_1\s_2\s_3^{-1}\s_4^{-1}\s_2\s_1^{-1}\s_2^{-1}
\quad\text{and}\quad
\xi_{\MM}=A_{2,4}A_{1,2}A_{3,4}A_{1,5}A_{3,5}.
\]

A {\it Macaulay~2} computation shows that $V_4(\MM)$ consists 
of sixteen $2$-torsion points.  
A {\it Mathematica} computation reveals that the Alexander polynomial 
of $\MM$ is an irreducible polynomial over $\Z$, consisting of $317$ 
monomials. The single variable Alexander polynomial is   
$\Delta_{\MM}(t)=(t-1)^5(t^2-t+1)(t^6-5t^5-t^4-6t^3-t^2-5t+1)$, 
and so $\delta(\MM)=0$.  Theorem~\ref{thm:2tors} gives 
$\Tors_{2,1}(\MM)=31$, the same as for $\LL$.  (Thus, 
the top variety of $\MM$ is not the union of translated subtori.)   
On the other hand, a computation yields $\Tors_{3,1}(\LL)=527$ and 
$\Tors_{3,1}(\MM)=421$, showing that $X(\LL)\not\simeq X(\MM)$.  


\section{Classification of $2$-Arrangements of $n\le 6$ Planes} 
\label{sec:classification}

We start by reviewing the rigid isotopy classification 
of arrangements of up to $6$ planes.  We then show that 
the invariants introduced in \S\ref{sec:detideals} are 
powerful enough to classify up to homotopy the complements of 
such arrangements.  

\subsection{}
As noted in \S\ref{rigidisotop}, rigid isotopy types of 
$2$-arrangements in $\R^4$ are in one-to-one correspondence with 
rigid isotopy types of skew-lines configurations in $\R^3$. 
An important concept introduced by Viro~\cite{V} was that of 
a mirror image of a configuration.  We now translate this notion 
to arrangements.

\begin{defn} 
An arrangement $\A'$ is called a {\em mirror image} of $\A$ 
if there is a reflection of $\R^4$ sending $\A$ to $\A'$.  
The mirror image of $\A$ is unique up to rigid isotopy; we 
denote it by $\overline{\A}$.  An arrangement $\A$ which is not isotopic 
to $\overline{\A}$ is called {\em non-mirror}.   
\end{defn}

As shown by Viro, there exist many non-mirror arrangements.  
For example, the complex arrangement, $\A_n$, and its mirror image 
under complex conjugation, $\overline{\A}_n$, are not rigidly isotopic 
provided $n\ge 3$.  Also, an arbitrary arrangement of $n$ lines is 
non-mirror, provided $n\equiv 3 \pmod{4}$.  

Viro~\cite{V} and Mazurovski\u{\i}~\cite{M1} classified, up to rigid 
isotopy, all configurations $6$ lines or less.  For up to $5$ lines,  
linking numbers invariants were used.  For $6$ lines, those invariants
cannot tell apart the $K$ and $L$ configurations.  For that, the Morton 
trace of the reduced full-braid is used in~\cite{M1}.  Translated to 
arrangements, the complete list of the $33$ rigid isotopy types is 
as follows:
\begin{align*}  \label{list6}
n=1: &\quad\A(1)\\
n=2: &\quad\A(12)\\
n=3: &\quad\A(123)^*\\
n=4: &\quad\A(1234)^*,\: \A(2134)\\
n=5: &\quad\A(12345)^*,\: \A(21345)^*,\: \A(21435)^*,\: \A(31425) \\
n=6: &\quad\A(123456)^*,\: \A(213456)^*,\: \A(321456),\: 
\A(214356)^*,\: \A(215436)^*, \\ 
     &\quad\A(312546),\:	\A(341256)^*,\: \A(314256)^*,\: 
\A(241536),\: \A(L)^*,\: \A(M)^*
\end{align*}
where $\A^*$ stands for a non-mirror arrangement $\A$ and its mirror 
image $\overline{\A}$.

\subsection{}
We now turn to the homotopy classification of complements of arrangements.  
Clearly, arrangements that are either rigidly isotopic, or mirror images 
of one another, have diffeomorphic (and thus, homotopy equivalent) complements.  
Thus, if we delete from the above list the mirror image $\overline{\A}$ from 
each pair $\A^*=(\A,\overline{\A})$, we are left with a list $20$ 
arrangements, such that, the complement of any arrangement of $n\le 6$ planes 
is homotopy equivalent to the complement of one of the arrangements in 
this shorter list.   Table~\ref{tab:classify} shows that there are no 
repetitions among the homotopy types of these $20$ arrangements.  
Hence, we have the following.  

\begin{thm}  \label{thm:classify6}
For $2$-arrangements of $n\le 6$ planes in $\R^4$, 
the homotopy types of complements are in one-to-one 
correspondence with the rigid isotopy types modulo mirror images. 
\end{thm}

\begin{table}
\begin{minipage}{1.02\columnwidth}
\[
\begin{array}{|c||c|c||c|c|c|}
\hline
\multicolumn{1}{|c||}{n}
&\multicolumn{1}{|c|}{\A}
&\multicolumn{1}{|c||}{\depth}
&\multicolumn{1}{|c|}{\Sigma_{n-2}}
&\multicolumn{1}{|c|}{\Tors_{2,1}}
&\multicolumn{1}{|c|}{\Tors_{3,1}}
\\ \hline\hline
1 & \A(1)        & 0 & 0  & 0 & 0 \\
\hline
2 & \A(12)       & 1 & 0  & 1 & 1 \\
\hline
3 & \A(123)      & 1 & 1 & 4 & 9\\
\hline
4 & \A(1234)     & 1 & 1  & 8 & 27 \\
& \A(2134)       & 2 & 2,3_2  & 8 & 45 \\
\hline
5 & \A(12345)    & 1 & 1 & 16 & 81 \\
& \A(21345)      & 2 & 2,3,4 & 16 & 135  \\
& \A(21435)      & 2 & 3,4_4  & 16 & 171 \\
& \A(31425)      & - & 5_{11} & 16 & 141 \\
\hline
6 & \A(123456)   & 1 & 1 & 32 & 243  \\
& \A(213456)     & 2 & 2,3,5  & 32 & 405 \\
& \A(321456)     & 2 & 2,4_2  & 32 & 405 \\
& \A(215436)     & 2 & 3,4_2,5_2  & 32 & 513 \\
& \A(214356)     & 2 & 3,4_2,5_3  & 32 & 513 \\
& \KK=\A(341256) & 3 & 4_3,5_3 & 32 & 567 \\ 
& \A(312546)     & 3 & 4,5_6,6 & 32 & 585 \\
& \A(314256)     & - & 5_6,6_5 & 32 & 495 \\
& \A(241536)     & - & 5_2,6_{13} & 31 & 513 \\ 
& \LL						      & - & 4_3,5_3 & 31 & 527 \\ 
& \MM 			     	  & - & 6_{16} & 31 & 421 \\
\hline
\end{array}
\]
\begin{center}
\caption{Arrangements of $n\le 6$ planes:  
Sequence of codimensions of components of $V_{n-2}$---where 
$i_k$ stands for $i$ repeated $k$ times---and number of 
$2$- and $3$-torsion points on $V_1$.\label{tab:classify}}
\vspace{-2pc}
\end{center}
\end{minipage}
\end{table}


\end{document}